\documentclass[12pt,a4paper]{amsart}
\usepackage{geometry}
\usepackage{latexsym}
\usepackage{bm}
\usepackage{amsmath}
\usepackage{amssymb}
\usepackage{amsmath,amscd}
\usepackage{tabls}
\usepackage{url}
\usepackage[utf8]{inputenc}
\usepackage{amsfonts}
\usepackage{amssymb}
\usepackage{mathrsfs}
\usepackage{amsmath}
\usepackage{amsthm}
\usepackage{amscd}
\usepackage{fancyhdr}
\usepackage{enumerate}
\usepackage{paralist}
\usepackage{xypic}
\usepackage{graphicx}
\usepackage{cite}
\usepackage{lipsum}
\usepackage{footmisc}
\usepackage[all,cmtip]{xy}

\theoremstyle{plain}
\newtheorem{Cor}{Corollary}[section]
\newtheorem{Def}[equation]{Definition}
\newtheorem{Thm}[equation]{Theorem}
\newtheorem{lem}[equation]{Lemma}
\newtheorem{prop}[equation]{Proposition}
\newtheorem{rem}[equation]{Remark}

\newtheorem{exa}[equation]{Example}

\numberwithin{equation}{section}

\title{Multiple Eisenstein L-Values}
\author{Zhongyu Jin}
\email{zyjin@pku.edu.cn}
\address{Zhongyu Jin \\School of Mathematical Sciences,
        Peking University,
         Beijing, China}

\begin{document}
\maketitle
\begin{abstract}
In this paper we study multiple Hecke L-values of Eisenstein series via the so-called iterated Eisenstein $\tau$-integrals and multiple Eisenstein $\tau$-series. As an application, we obtain an explicit relationship between double Hecke L-values of Eisenstein series and holomorphic double modular values of Eisenstein series defined by Brown.
\end{abstract}

\section{Introduction}

Manin \cite{Y.M1} constructed iterated integrals of cusp forms. Furthermore, Brown \cite{F.B2} generalized Manin's work to general modular forms. Choie and Ihara \cite{K.I1} studied multiple Hecke L-functions of cusp forms, which are connected to Manin's work. In this paper, we study the multiple Eisenstein L-values, which are the values of multiple Hecke L-functions of Eisenstein series at positive integers, with the help of Brown's work.

\subsection{Iterated integrals of cusp forms}

Suppose that $f_{i}(\tau)=\sum_{n\geq1}a_{n}(f_{i})e^{2\pi i\tau}$ are cusp forms of weight $2k_{i}$ for the modular group $SL_{2}(\mathbb{Z})$, $i=1,\cdots,r$. Manin \cite{Y.M1}, \cite{Y.M2} studied the following function on the upper half complex plane:
$$J(\tau)=\int_{\tau}^{i\infty}f_{1}(\tau_{1})(X-Y\tau_{1})^{2k_{1}-2}\cdots f_{r}(\tau_{r})(X-Y\tau_{r})^{2k_{r}-2}d\tau_{r}\cdots d\tau_{1},$$
where $X,Y$ are variables and the integral is an iterated integral. It does not depend on the path from $\tau$ to $i\infty$ and thus well-defined. The function $J(\tau)$ is a homogeneous polynomial of $X$ and $Y$, whose coefficients are linear combinations of
$$\sum_{n_{1},\cdots,n_{r}>0}\frac{a_{n_{1}}(f_{1})a_{n_{2}}(f_{2})\cdots a_{n_{r}}(f_{r})}{(n_{1}+\cdots+n_{r})^{\alpha_{1}}(n_{2}+\cdots+n_{r})^{\alpha_{2}}\cdots n_{r}^{\alpha_{r}}}e^{2\pi i(n_{1}+\cdots n_{r})}$$
for positive integers $\alpha_{i}\in\mathbb{Z}$, $i=1,\cdots,r$. When $\tau=0$, $J(0)$ is also well-defined and we call it the iterated Eichler-Shimura integral of the cusp forms $f_{i}(\tau)$, $i=1,\cdots,r$. This is a homogeneous polynomial of $X, Y$ and can be regarded as the multiple period polynomial of cusp forms. The coefficients of the polynomial $J(0)$ are linear combinations of
$$\sum_{n_{1},\cdots,n_{r}>0}\frac{a_{n_{1}}(f_{1})a_{n_{2}}(f_{2})\cdots a_{n_{r}}(f_{r})}{(n_{1}+\cdots+n_{r})^{\alpha_{1}}(n_{2}+\cdots+n_{r})^{\alpha_{2}}\cdots n_{r}^{\alpha_{r}}}.$$
Such a series is the multiple Hecke L-value of the cusp forms $f_{i}(\tau)$ (up to a power of $2\pi i$), $i=1,\cdots,r$, at integer points $\alpha_{1},\cdots,\alpha_{r}$ as in \cite{K.I1}. Manin \cite{Y.M1} studied their shuffle relations.

\subsection{Eisenstein series}

It is natural for us to consider not only cusp forms, but also Eisenstein series. These are the objects in this paper. Let
$$E_{2k}(\tau)=-\frac{b_{2k}}{4k}+\sum_{n\geq1}\sigma_{2k-1}(n)q^{n}$$
be the (Hecke normalized) Eisenstein series of weight $2k$ for $SL_{2}(\mathbb{Z})$, where $k\geq 2$ is a positive integer, $b_{2k}$ is the $2k$-th Bernoulli number and $\sigma_{2k-1}$ is the $(2k-1)$-th divisor function. Denote by
$$E_{2k}^{0}(\tau)=\sum_{n\geq1}\sigma_{2k-1}(n)q^{n}$$
the cuspital part of the Eisenstein series. Consider the so-called iterated Eisenstein $\tau$-integral as
$$Int(E_{2k_{1}}^{0},\cdots,E_{2k_{r}}^{0};\alpha_{1},\cdots,\alpha_{r})(\tau)$$
$$=\int_{\tau<\tau_{1}<\cdots<\tau_{n}<i\infty}E_{2k_{1}}^{0}(\tau_{1})\tau_{1}^{\alpha_{1}-1}\cdots E_{2k_{r}}^{0}(\tau_{r})\tau_{r}^{\alpha_{r}-1}d\tau_{r}\cdots d\tau_{1}.$$
This iterated integral is well-defined since $\tau$ is in the upper half complex plane and $E_{2k}^{0}(\tau)$ has no singular part at $i\infty$ for $k\geq2$. We denote the space of all iterated Eisenstein $\tau$-integrals by $\mathcal{IEI}_{\tau}$.

We also consider the following function of $\tau$ on the upper half complex plane as
$$L^{(t)}(E_{2k_{1}}^{0},\cdots,E_{2k_{r}}^{0};\alpha_{1},\cdots,\alpha_{r})(\tau)$$
$$=\tau^{t}(2\pi i)^{-\alpha_{1}-\cdots-\alpha_{r}}\sum\limits_{n_{1},\cdots,n_{r}>0}\frac{\sigma_{2k_{1}-1}(n_{1})\sigma_{2k_{2}-1}(n_{2})\cdots\sigma_{2k_{r}-1}(n_{r})}
{(n_{1}+\cdots+n_{r})^{\alpha_{1}}(n_{2}+\cdots+n_{r})^{\alpha_{2}}\cdots n_{r}^{\alpha_{r}}}e^{2\pi i(n_{1}+\cdots+n_{r})\tau}$$
for positive integers $\alpha_{1},\cdots,\alpha_{r}$, which we call a multiple Eisenstein $\tau$-series. It is also well-defined. We add the non-negative integer $t$ in the definition in order to state the results below simpler. Denote by $\mathcal{MEL}_{\tau}$ the space of all multiple Eisenstein $\tau$-series.

By induction on $r$ and using the operators $(1+\frac{\partial}{\partial \tau})^{-1}$, we express multiple Eisenstein $\tau$-series as linear combinations of iterated Eisenstein $\tau$-integrals. Also we give the converse expression using the operator $\frac{\partial}{\partial \tau}$. In particular, we prove the following theorem.

\begin{Thm}\label{M1}
The $\mathbb{Q}$-vector space $\mathcal{MEL}_{\tau}$ of all multiple Eisenstein $\tau$-series is a $\mathbb{Q}$-algebra with the length filtration $L^{\bullet}$ under the natural product, and we have
$$L^{n_{1}}(\mathcal{MEL}_{\tau})\cdot L^{n_{2}}(\mathcal{MEL}_{\tau})\subset L^{n_{1}+n_{2}}(\mathcal{MEL}_{\tau}).$$
Also there is an algebra equation $\mathcal{MEL}_{\tau}=\mathcal{IEI}_{\tau}$.
\end{Thm}

A result in \cite{M.D} gives a way to prove the linear independence of special types of functions, we will recall it in Section $3$. As an application, authors of \cite{N.M} proved the linear independence of their regularized iterated Eiseistein integrals by considering the Fourier coefficients of Eisenstein series (note that their integrals are similar, but not totally the same as the iterated Eisenstein $\tau$-integrals here). As an analogue, we get a basis for the $\mathbb{Q}$-vector space $\mathcal{MEL}_{\tau}$, and thus the $\mathbb{Q}$-vector space of iterated Eisenstein $\tau$-integrals. Precisely, we have:

\begin{Thm}\label{M2}
As functions on the upper half complex plane, the set of multiple Eisenstein $\tau$-series
$$\{L^{(t)}(E_{2k_{1}}^{0},\cdots,E_{2k_{r}}^{0};\alpha_{1},\cdots,\alpha_{r})(\tau);r\geq0,k_{i}\geq2,\alpha_{i}\geq1\}$$
are $\mathbb{C}$-linear independent functions. As a consequence, this set forms a basis of the $\mathbb{Q}$-vector space $\mathcal{MEL}_{\tau}$.
\end{Thm}

\subsection{Multiple Eisenstein L-values}

It is obvious to see that when $\tau\rightarrow 0$ and $t=0$, the multiple Eisenstein L-series
$$L^{(t)}(E_{2k_{1}}^{0},\cdots,E_{2k_{r}}^{0};\alpha_{1},\cdots,\alpha_{r})(\tau)$$
formally reduces to
$$L(E_{2k_{1}}^{0},\cdots,E_{2k_{r}}^{0};\alpha_{1},\cdots,\alpha_{r})$$
$$=(2\pi i)^{-\alpha_{1}-\cdots-\alpha_{r}}\sum\limits_{n_{1},\cdots,n_{r}>0}\frac{\sigma_{2k_{1}-1}(n_{1})\sigma_{2k_{2}-1}(n_{2})\cdots\sigma_{2k_{r}-1}(n_{r})}
{(n_{1}+\cdots+n_{r})^{\alpha_{1}}(n_{2}+\cdots+n_{r})^{\alpha_{2}}\cdots n_{r}^{\alpha_{r}}}.$$
In particular, when $r=1$, we have $L(E_{2k}^{0};\alpha)=(2\pi i)^{-\alpha}\zeta(\alpha)\zeta(\alpha-2k+1)$ for $\alpha\neq 0$ and $\alpha\neq 2k$, where $\zeta(s)$ is the Riemman zeta function, and the meromorphic extension $L(E_{2k}^{0};s)$ has simple poles at $\alpha=0$ and $\alpha=2k$.

The series above is not directly well-defined. However, we have the following result. For $k=1,\cdots,r$, suppose that $\phi_{k}(s)=\sum_{n=1}^{\infty}\frac{a_{k}(n)}{n^{s}}$ is convergent absolutely for $Re(s)>\alpha_{k}>0$ and can be continued meromorphically to the whole complex plane, holomorphic except for a possible pole at $s=\alpha_{k}$ of order at most $1$. Besides, assume that $\phi_{k}(s)=O(|Im(s)|^{A})$ as $|s|\rightarrow \infty$ for some non-negative constant $A$. Then the following theorem holds.

\begin{Thm}[Matsumoto and Tanigawa \cite{K.M}]
The multiple series
$$\sum_{m_{1},\cdots,m_{r}=1}^{\infty}\frac{a_{1}(m_{1})a_{2}(m_{r})\cdots a_{r}(m_{r})}{(m_{1}+\cdots+m_{r})^{s_{1}}(m_{2}+\cdots+m_{r})^{s_{2}}\cdots m_{r}^{s_{r}}}$$
can be continued meromorphically to the whole $\mathbb{C}^{r}$ space, and its possible singularities are located only on the subsets $\mathbb{C}^{r}$, each of which is defined by one of the following equations:
$$s_{j}+\cdots+s_{r}=\alpha_{j}+\delta_{j+1}\alpha_{j+1}+\cdots+\delta_{r}\alpha_{r}-n,$$
where $1\leq j\leq r$, $n\in\mathbb{Z}$ non-negative and $\delta_{k}=0$ or $1$ for $2\leq k\leq r$.
Furthermore, if $\phi_{k}(s)$ are entire on $\mathbb{C}$ for $k=1,\cdots,r$, then the multiple series is holomorphic with respect to $s_{1},\cdots,s_{r}$.
\end{Thm}

When $\phi_{k}$ are cusp forms for $k=1,\cdots,r$, we call the series multiple Hecke L-function as in \cite{K.I1}, it is holomorphic for $s_{1},\cdots,s_{r}$ and its values at integer points occur in Manin's work. For Eisenstein series,
$$L(E_{2k_{1}}^{0},\cdots,E_{2k_{r}}^{0};\alpha_{1},\cdots,\alpha_{r})$$
is the value of the multiple Hecke L-functions of Eisenstein series at integral points $\alpha_{i}\in\mathbb{Z}_{+}$, and at the singular points, we can regard it as a Laurent series of $\alpha_{i}$. By abusing the name, we call it a multiple Eisenstein L-value for convenience. As a direct corollary of Theorem \ref{M1}, we have:

\begin{Cor}
The $\mathbb{Q}$-vector space $\mathcal{MEL}$ of multiple Eisenstein L-values is a $\mathbb{Q}$-algebra with a length filtration $L^{\bullet}$ satisfying that
$$L^{n_{1}}(\mathcal{MEL})\cdot L^{n_{2}}(\mathcal{MEL})\subset L^{n_{1}+n_{2}}(\mathcal{MEL}).$$
Furthermore, $Int(E_{2k_{1}}^{0},\cdots,E_{2k_{r}}^{0};\alpha_{1},\cdots,\alpha_{r})(0)$ has a Laurent expression after the meromorphic extension.
\end{Cor}

\subsection{Holomorphic multiple modular values}

Brown \cite{F.B2} defined iterated Eichler integral generalizing Manin's work. By this, he also defined (homomorphic) multiple modular values of modular forms for the group $SL_{2}(\mathbb{Z})$. When we only consider cusp forms, Brown's construction is compatible with Manin's. As for Eisenstein series, we can find interesting information among them, including Riemann zeta values and L-values of cusp forms outside the critical line. More details will be introduces in Section $2$.

The relationship between multiple Hecke L-values and iterated integrals of cusp forms has been studied by Choie and Ihara \cite{K.I1}, we give the statement in Section $3$. Here we want to set up the relationship between multiple Eisenstein L-values and holomorphic multiple modular values of Eisenstein series. In order to do this, we need to use the previous results, and in length two, we have:

\begin{Thm}\label{Intro 3}
A holomorphic double modular value of Eisenstein series can be expressed as a $\mathbb{Q}[2\pi i]$-linear combination of double Eisenstein $L$-values explicitly.
\end{Thm}

It has been shown that $L(E_{2k}^{0};\alpha)$ is a product of two Riemman zeta values outside $\alpha=0$ and $\alpha=2k$. In the multiple case, multiple Eisenstein L-values are connected to multiple zeta values, one can find more relative work in \cite{F.B2}, \cite{B.E} and \cite{N.M}. Then Theorem \ref{Intro 3} shows an explicit relationship between double modular values and double Eisenstein L-values, and thus double zeta values.

\section{Regularized Iterated Integrals and Holomorphic Multiple Modular Values}

In this section, we recall the regularized iterated integrals of modular forms for $SL_{2}(\mathbb{Z})$ and (holomorphic) multiple modular valued defined by Brown \cite{F.B2}. One can also find relevant results in \cite{Y.M1} and \cite{Y.M2} of cusp forms.

\subsection{Regularized iterated integrals}

Denote by $\mathfrak{h}$ the upper half complex plane and by $\Gamma=SL_{2}(\mathbb{Z})$, which is the modular group generated by
$$T=\begin{pmatrix}
1 & 1 \\
0 & 1
\end{pmatrix}, S=\begin{pmatrix}
0 & -1 \\
1 & 0
\end{pmatrix}.$$
Let $\mathcal{M}_{k}(\Gamma)$ be the space of modular forms for $\Gamma$ of weight $k$. Fix a rational basis $\mathcal{B}$ of $\mathcal{M}(\Gamma)=\bigoplus_{k}\mathcal{M}_{k}(\Gamma)$. We assume that $\mathcal{B}=\cup_{k}\mathcal{B}_{k}$, where $\mathcal{B}_{k}$ is a basis of $\mathcal{M}_{k}(\Gamma)$, and that $\mathcal{B}_{k}$ is compatible with the action of Hecke operators. We always suppose that $\mathcal{B}_{2k}$ contains the Hecke normalised Eisenstein series $E_{2k}(\tau)$.

Define a $\mathbb{Q}$-vector space with a basis consisting of certain symbols indexed by $\mathcal{B}_{k}$ to be
$$M_{k}=\langle a_{f},f\in\mathcal{B}_{k}\rangle_{\mathbb{Q}}$$
for any $k\geq2$. Their dual spaces are denoted by
$$M_{k}^{\vee}=\langle A_{f},f\in\mathcal{B}_{k}\rangle_{\mathbb{Q}},$$
where $\langle a_{f}, A_{g}\rangle=\delta_{f,g}$ and $\delta$ is the Kronecker delta function.

Denote by $V_{k}$ the space of homogeneous polynomials of variables $X$ and $Y$ of degree $k$, there is a natural $\Gamma$ action on it as
$$\gamma(X,Y)=(aX+bY,cX+dY), \forall\gamma=\begin{pmatrix}
a & b \\
c & d
\end{pmatrix}\in\Gamma.$$
Then we have the graded right $SL_{2}(\mathbb{Z})$-module
$$M^{\vee}=\bigoplus_{k\geq2}M_{k}^{\vee}\otimes V_{k-2},$$
which has one copy of $V_{k-2}$ for every element of $\mathcal{B}_{k}$, and the group $\Gamma$ acts on $M_{k}^{\vee}$ trivially. Define
$$\mathcal{U}_{1,1}^{dR,hol}(\mathbb{C})=\{S\in \mathbb{C}\langle\langle M^{\vee}\rangle\rangle^{\times}, S\ is\ group\ like\},$$
where $\mathbb{C}\langle\langle M^{\vee}\rangle\rangle$ is the ring of formal power series of elements in $M^{\vee}$ with complex coefficients.

For any modular form $f(\tau)\in \mathcal{M}_{2k}(\Gamma)$, define the differential forms to be
$$\underline{f}(\tau)=(2\pi i)^{2k-1}f(\tau)(X-Y\tau)^{2k-2}d\tau,$$
$$\underline{f}^{\infty}(\tau)=(2\pi i)^{2k-1}a_{f}(0)(X-Y\tau)^{2k-2}d\tau,$$
where $a_{f}(0)$ is the Fourier constant term of $f(\tau)$. Then for any $\gamma\in\Gamma$, we have $\underline{f}(\gamma(\tau))|_{\gamma}=\underline{f}(\tau)$, where $|_{\gamma}$ means the action of $\Gamma$ on $V_{2k-2}$ as above. Moreover, we can define differential forms
$$\Omega(\tau)=\sum_{f\in\mathcal{B}}A_{f}\underline{f}(\tau), \Omega^{\infty}(\tau)=\sum_{f\in\mathcal{B}}A_{f}\underline{f}^{\infty}(\tau).$$
We also have $\Omega(\gamma(\tau))|_{\gamma}=\Omega(\tau)$ for any $\gamma\in\Gamma$. For any two points $\tau_{1},\tau_{2}\in\mathfrak{h}$, consider the iterated integral
$$I(\tau_{1},\tau_{2})=1+\int_{\tau_{1}}^{\tau_{2}}\Omega(\tau)+\int_{\tau_{1}}^{\tau_{2}}\Omega(\tau)\Omega(\tau)+\cdots,$$
$$I^{\infty}(\tau_{1},\tau_{2})=1+\int_{\tau_{1}}^{\tau_{2}}\Omega^{\infty}(\tau)+\int_{\tau_{1}}^{\tau_{2}}\Omega^{\infty}(\tau)\Omega^{\infty}(\tau)+\cdots.$$
Since $\Omega(\tau)$ and $\Omega^{\infty}(\tau)$ are integral and the upper half complex plane $\mathfrak{h}$ is simply connected, the iterated integrals above do not depend on the path from $\tau_{1}$ to $\tau_{2}$. Brown \cite{F.B2} proved the following result and give the definition as:

\begin{Def}
For any $\tau\in\mathfrak{h}$, the limit below is well-defined and we call it the iterated Eichler integral from $\tau$ to $i\infty$:
$$I(\tau,i\infty)=\lim_{z\rightarrow i\infty}I(\tau,z)I^{\infty}(z,0).$$
\end{Def}

We may also call $I(\tau,i\infty)$ the regularized iterated integral of modular forms from $\tau$ to $i\infty$. If we replace $i\infty$ by any point $\alpha\in\mathbb{Q}$, we can still define the integral. In this case, take $\gamma\in\Gamma$ such that $\gamma(\alpha)=i\infty$ and define
$$I(\tau,\alpha)=I(\gamma(\tau),\gamma(\alpha))|_{\gamma}.$$
Also we can define
$$I(\alpha,i\infty)=I(\alpha,\tau)I(\tau,i\infty),$$
by choosing any $\tau\in\mathfrak{h}$. It does not depend on the choice of $\tau$ and thus well-defined by the following fundamental proposition:

\begin{prop}
The regularized iterated integrals satisfy the following properties:

$(1)$. Differential property: $\partial_{\tau}I(\tau,i\infty)=-\Omega(\tau)I(\tau,i\infty)$.

$(2)$. For any $\tau_{1},\tau_{2},\tau_{3}\in\mathfrak{h}\cup\mathbb{Q}\cup\{i\infty\}$, we have $I(\tau_{1},\tau_{3})=I(\tau_{1},\tau_{2})I(\tau_{2},\tau_{3})$.

$(3)$. Modular property: for any $\gamma\in\Gamma$, we have $I(\gamma(0),\gamma(i\infty))|_{\gamma}=I(0,i\infty)$.

$(4)$. Shuffle property: the integral $I(\tau,i\infty)$ is shuffled and invertible.
\end{prop}

\subsection{Holomorphic multiple modular values}

With the help of the iterated Eichler integral, Brown \cite{F.B2} proved the following result and defined holomorphic multiple modular values.

\begin{prop}[Definition]
For every $\gamma\in\Gamma$, there exists a series $\mathcal{C}_{\gamma}\in\mathcal{U}^{dR,hol}_{1,1}(\mathbb{C})$ such that
$$I(\tau;i\infty)=I(\gamma(\tau),i\infty)|_{\gamma}\mathcal{C}_{\gamma}$$
It does not depend on $\tau$. It satisfies the cocycle relation
$$\mathcal{C}_{gh}=\mathcal{C}_{g}|_{h}\mathcal{C}_{h}, \forall g,h\in\Gamma.$$
\end{prop}

\begin{rem}
Actually, $\mathcal{C}_{\gamma}=I(\gamma_{-1}(i\infty),i\infty)$.
\end{rem}

\begin{Def}
Define the ring of holomorphic multiple modular values $\mathcal{MMV}^{hol}_{\Gamma}$ for $\Gamma$ to be the $\mathbb{Q}$-algebra generated by the coefficients of
$$A_{f_{1}}\cdots A_{f_{n}}X_{1}^{i_{1}-1}Y_{1}^{2k_{1}-i_{1}-1}\cdots X_{n}^{i_{n}-1}Y_{n}^{2k_{n}-i_{n}-1}$$ for any non-negative $n\in\mathbb{Z}$ and $f_{i}\in\mathcal{B}_{2k_{i}}$ in $\mathcal{C}_{\gamma}$, $\forall\gamma\in\Gamma$.
\end{Def}

\begin{rem}
We call elements in $\mathcal{MMV}_{\Gamma}^{hol}$ holomorphic to distinguish from the algebraic ones constructed in \cite{F.B2}. Since we only consider the homomorphic ones here, we would like to call them multiple modular values and denote $\mathcal{MMV}_{\Gamma}^{hol}$ by $\mathcal{MMV}$ for convenience.
\end{rem}

Since the group $\Gamma$ is generated by $S$ and $T$, the ring $\mathcal{MMV}$ only depends on the coefficients of $\mathcal{C}_{S}$ and $\mathcal{C}_{T}$ because of the cocycle properties of $\mathcal{C}$. As in \cite{F.B2}, we can calculate $\mathcal{C}_{T}$ directly by
$$\mathcal{C}_{T}=I^{\infty}(-1,0)$$
and thus the coefficients of $\mathcal{C}_{T}$ belong to $\mathbb{Q}[2\pi i]$. They can be regarded as periods of the motivic fundamental group of $\mathbb{G}_{m}$. As for $\mathcal{C}_{S}$, take $\tau=i$ and $\gamma=S$, we have
$$\mathcal{C}_{S}=I(i,i\infty)|^{-1}_{S}I(i,i\infty).$$
There are much interesting information hidden in it.

The computation of length one parts of $\mathcal{C}_{S}$ is classical, Brown \cite{F.B2} gave a set of rational cocycles to describe them. Let $c(x)=\frac{1}{e^{x}-1}$, define a set of rational cocycles $e_{2k}\in Z^{1}(\Gamma;V_{2k-2})$ via the generation series:
$$e^{0}=\sum_{k\geq2}\frac{2}{(2k-2)!}e^{0}_{2k},$$
where $e^{0}$ is the unique cocycle in $V_{\infty}$ defined on $\Gamma$ by
$$e^{0}(S)=c(X)c(Y),\ e^{0}(T)=\frac{1}{Y}(c(X+Y)-c(X))-\frac{1}{12}.$$
Notice that these $e^{0}_{2k}$ satisfy the abelian cocycle relations, and explicitly, for $k\geq2$, we have:
$$e^{0}_{2k}(S)=\frac{(2k-2)!}{2}\sum_{i=1}^{k-1}\frac{b_{2i}}{(2i)!}\frac{b_{2k-2i}}{(2k-2i)!}X^{2i-1}Y^{2k-2i-1},$$
$$e^{0}_{2k}(T)=\frac{(2k-2)!}{2}\frac{b_{2k}}{(2k)!}\frac{(X+Y)^{2k-1}-X^{2k-1}}{Y}.$$

\begin{Thm}[Brown \cite{F.B2}]
The following formulas hold for $\mathcal{C}_{S}$:

$(1)$. For any cusp form $f$ of weight $k$,
$$\mathcal{C}_{S}(f)=(2\pi i)^{k-1}\int_{0}^{i\infty}f(\tau)(X-\tau Y)^{k-2}d\tau.$$
Thus it is the period polynomial of $f$.

$(2)$. For Hecke normalized Eisenstein series $E_{2k}$, we have:
$$\mathcal{C}_{S}(E_{2k})=(2\pi i)^{2k-1}e^{0}_{2k}(S)+\frac{(2k-2)!}{2}\zeta(2k-1)(X^{2k-2}-Y^{2k-2}).$$
\end{Thm}

As for multiple modular values of length $2$, the story becomes more complex. Let $a,b \geq0$ and $k\geq0$, define
$$I^{2k}_{2a,2b}=\partial^{k}Im(\mathcal{C}_{e_{2a}e_{2b}}+\overline{b}_{2a}\cup\overline{e}^{0}_{2b}-\overline{e}^{0}_{2a}\cup\overline{b}_{2b}),$$
where $\mathcal{C}_{e_{2a}e_{2b}}$ is the coefficient of $A_{E_{2a}}A_{E_{2b}}$ in $\mathcal{C}$, it defines a map from $\Gamma$ to the space of homogeneous polynomials of $X, Y$. The notation $Im$ means the imaginary part, $\cup$ means the cup product of cocycles and
$$\overline{b}_{2k}=\frac{(2k-2)!}{2}\zeta(2k-1)Y^{2k-1};\ \overline{e}^{0}_{2k}=(2\pi i)^{2k-1}e^{0}_{2k},$$
$$\partial^{k}:\mathbb{Q}[X_{1},X_{2},Y_{1},Y_{2}]\rightarrow \mathbb{Q}[X,Y],\ \partial^{k}=\pi_{d}\circ(\partial_{12}^{k}),$$
here
$$\pi_{d}:\mathbb{Q}[X_{1},X_{2},Y_{1},Y_{2}]\rightarrow\mathbb{Q}[X,Y];(X_{i},Y_{i})\rightarrow (X,Y), i=1,2$$
$$\partial_{12}=\frac{\partial}{\partial_{X_{1}}}\frac{\partial}{\partial_{Y_{2}}}-\frac{\partial}{\partial_{X_{2}}}\frac{\partial}{\partial_{Y_{1}}}.$$
Then the following theorem holds.

\begin{Thm}[Brown \cite{F.B2}]\label{Brown's l2 calculation}
Fix $a,b,k$ as above, let $w=2a+2b-2k-2$, the cochain $I^{2k}_{2a,2b}$ is a cocycle, it means $I^{2k}_{2a,2b}\in Z^{1}(\Gamma, V_{2a+2b-2k-4})$, and we have
$$I^{2k}_{2a,2b}(S)\equiv \sum_{\{g\}}(2\pi i)^{k}\Lambda(g,w+k)P_{g}^{\epsilon}(X,Y)\ (mod\ \delta^{0}(V_{w-2}\otimes \mathbb{C})_{S}),$$
where the sum ranges over a basis of Hecke normalised cusp eigenforms of weight $w$, $P^{\varepsilon}_{g}\in P_{w-2}\otimes K_{g}$ are Hecke-invariance period polynomials, $\varepsilon\in\{+,-\}$ denotes the conjugate-invariance if $k$ is odd and conjugate-anti-invariant if $k$ is even, and $K_{g}$ is the field generated by the Fourier coefficients of $g$.
\end{Thm}

An interesting fact following from Brown's theorem is that some L-values of Hecke normalized cusp forms outside the critical line appear in $\mathcal{MMV}$. We will give some direct calculation and show their relationship among multiple Eisenstein L-values of length $2$ in Section $4$.

\section{Iterated Eisenstein $\tau$-integrals and Multiple Eisenstein $\tau$-series}

In this section, we define and study iterated Eisenstein $\tau$-integrals and multiple Eisenstein $\tau$-series. We will prove Theorem \ref{M1} and Theorem \ref{M2} stated in the introduction. The definitions are the analogues of ones in \cite{Y.M1} and the strategies come from \cite{K.I1}, \cite{N.M} and \cite{Y.M1}.

\subsection{Definitions}

For any positive integer $k\geq2$ and the Hecke normalized Eisenstein series $E_{2k}(\tau)$, denote by
$$E^{0}_{2k}(\tau)=\sum_{n\geq0}\sigma_{2k-1}(n)q^{n}\ \text{and}\ E^{\infty}_{2k}(\tau)=-\frac{b_{2k}}{4k}$$
the cuspital part and Fourier constant term of $E_{2k}(\tau)$ respectively. For $\tau\in\mathfrak{h}$ and integers $r\geq0$, $k_{1},\cdots,k_{r}\geq2$, $\alpha_{1},\cdots,\alpha_{r}\geq1$, we define the iterated Eisenstein $\tau$-integral as
$$Int(E_{2k_{1}}^{0},\cdots,E_{2k_{r}}^{0};\alpha_{1},\cdots,\alpha_{r})(\tau)$$
$$=\int_{\tau<\tau_{1}<\cdots<\tau_{n}<i\infty}E_{2k_{1}}^{0}(\tau_{1})\tau_{1}^{\alpha_{1}-1}\cdots E_{2k_{r}}^{0}(\tau_{r})\tau_{r}^{\alpha_{r}-1}d\tau_{r}\cdots d\tau_{1}.$$
As stated in the introduction, this iterated integral is well-defined. We call $r$ the length of it, and when $r=0$, we let the integral equal to $1$.

\begin{rem}
Lochak, Matthes and Schneps \cite{N.M} studied the regularized iterated integrals of Eisenstein series. The weight $0$ Eisenstein series, i.e., $E_{0}(\tau)=-1$, is also considered in their paper. Actually, we can see that the iterated Eisenstein $\tau$-integrals can be written as $\mathbb{Q}$-linear combinations of regularized iterated integrals of Eisenstein series as in \cite{N.M}.
\end{rem}

\begin{rem}
The iterated Eisenstein $\tau$-integrals are closely connected to elliptic multiple zeta values, more detail can be found in \cite{N.M}, \cite{N.M2} and so on.
\end{rem}

\begin{Def}
For $\tau \in \mathfrak{h}$, denote by $\mathcal{IEI_{\tau}}$ the $\mathbb{Q}$-vector space
$$\mathcal{IEI_{\tau}}=\langle Int(E_{2k_{1}}^{0},\cdots,E_{2k_{r}}^{0};\alpha_{1},\cdots,\alpha_{r})(\tau); k_{i}\geq2,\alpha_{i}\geq1, r\geq0\rangle,$$
where $\langle\alpha,\beta,\cdots\rangle$ denotes the $\mathbb{Q}$-vector space spanned by $\alpha, \beta,\cdots$. The length filtration of $\mathcal{IEI}_{\tau}$ is defined as follows:
$$L^{n}(\mathcal{IEI_{\tau}})=\langle Int(E_{2k_{1}}^{0},\cdots,E_{2k_{r}}^{0};\alpha_{1},\cdots,\alpha_{r})(\tau) \in \mathcal{IEI(\tau)}; r\leq n\rangle.$$
\end{Def}

By the shuffle product among iterated integrals, $\mathcal{IEI_{\tau}}$ is a $\mathbb{Q}$-algebra and it is routine to see
$$L^{n_{1}}(\mathcal{IEI}_{\tau})\cdot L^{n_{2}}(\mathcal{IEI}_{\tau})\subset L^{n_{1}+n_{2}}(\mathcal{IEI}_{\tau}).$$
There is also a natural operator $\frac{\partial}{\partial \tau}=(2\pi i)q\frac{\partial}{\partial q}$ acting on these iterated Eisenstein $\tau$-integrals as
$$\frac{\partial}{\partial\tau}Int(E_{2k_{1}}^{0},\cdots,E_{2k_{r}}^{0};\alpha_{1},\cdots,\alpha_{r})(\tau)$$
$$=-E_{2k_{1}}^{0}(\tau)\tau^{\alpha_{1}-1}Int(E_{2k_{2}}^{0},\cdots,E_{2k_{r}}^{0};\alpha_{2},\cdots,\alpha_{r})(\tau),$$
where we let $Int(E_{2k_{2}}^{0},\cdots,E_{2k_{r}}^{0};\alpha_{2},\cdots,\alpha_{r})(\tau)=1$ when $r=1$.

Next, we consider multiple Eisenstein $\tau$-series. For $\tau\in \mathfrak{h}$ and integers $r\geq0$, $k_{1},\cdots,k_{r}\geq2$, $\alpha_{1},\cdots,\alpha_{r}\geq1$, $t\geq0$, define the multiple Eisenstein $\tau$-series to be the series
$$L^{(t)}(E_{2k_{1}}^{0},\cdots,E_{2k_{r}}^{0};\alpha_{1},\cdots,\alpha_{r})(\tau)$$
$$=\tau^{t}(2\pi i)^{-\alpha_{1}-\cdots-\alpha_{r}}\sum\limits_{n_{1},\cdots,n_{r}>0}\frac{\sigma_{2k_{1}-1}(n_{1})\sigma_{2k_{2}-1}(n_{2})\cdots\sigma_{2k_{r}-1}(n_{r})}
{(n_{1}+\cdots+n_{r})^{\alpha_{1}}(n_{2}+\cdots+n_{r})^{\alpha_{2}}\cdots n_{r}^{\alpha_{r}}}e^{2\pi i(n_{1}+\cdots+n_{r})\tau},$$
and when $r=0$ we let it equal to $1$. Since $Im(\tau)>0$ and $\alpha_{i}\geq1$, this series is convergent and thus well-defined. We call $r$ the length of this series.

The natural operator $\frac{\partial}{\partial \tau}=(2\pi i)q\frac{\partial}{\partial_{q}}$ also acts on the multiple Eisenstein $\tau$-series and we have
$$\frac{\partial}{\partial\tau}L^{(t)}(E_{2k_{1}}^{0},\cdots,E_{2k_{r}}^{0};\alpha_{1},\cdots,\alpha_{r})(\tau)$$
$$=tL^{(t-1)}(E_{2k_{1}}^{0},\cdots,E_{2k_{r}}^{0};\alpha_{1},\cdots,\alpha_{r})(\tau)
+L^{(t)}(E_{2k_{1}}^{0},\cdots,E_{2k_{r}}^{0};\alpha_{1}-1,\cdots,\alpha_{r})(\tau)$$
for $\alpha_{1}\geq2$. Here we let $L^{(t-1)}(E_{2k_{1}}^{0},\cdots,E_{2k_{r}}^{0};\alpha_{1},\cdots,\alpha_{r})(\tau)=0$ when $t=0$.

\begin{Def}
Denote by $\mathcal{MEL_{\tau}}$ the $\mathbb{Q}$-vector space
$$\mathcal{MEL_{\tau}}=\langle L^{(t)}(E_{2k_{1}}^{0},\cdots,E_{2k_{r}}^{0};\alpha_{1},\cdots,\alpha_{r})(\tau);k_{i}\geq2, \alpha_{i}\geq1, r\geq0, t\geq0 \rangle.$$
The length filtration of it is defined as:
$$L^{n}(\mathcal{MEL_{\tau}})=\langle L^{(t)}(E_{2k_{1}}^{0},\cdots,E_{2k_{r}}^{0};\alpha_{1},\cdots,\alpha_{r})(\tau) \in \mathcal{IEI(\tau)}; r\leq n\rangle.$$
\end{Def}

For any positive integers $i,j\geq1$, $i+j=k$ and $m,n\geq1$, the following formula is well-known and will be used again below:
$$\frac{1}{m^{i}n^{j}}=\sum_{r+s=k}\frac{{r-1 \choose i-1}}{(m+n)^{r}n^{s}}+\frac{{r-1 \choose j-1}}{(m+n)^{r}m^{s}}.$$
With the help of this formula, we give an example of the relation among multiple Eisenstein $\tau$-series of length $2$.

\begin{exa}\label{Stuffle relation among DEL}
Note that
$$L^{(0)}(E_{2k_{1}}^{0},E_{2k_{2}}^{0};\alpha_{1},\alpha_{2})(\tau)
=(2\pi i)^{-\alpha_{1}-\alpha_{2}}
\sum_{n_{1},n_{2}>0}\frac{\sigma_{2k_{1}-1}(n_{1})\sigma_{2k_{2}-1}(n_{2})}{(n_{1}+n_{2})^{\alpha_{1}}n_{2}^{\alpha_{2}}}q^{n_{1}+n_{2}},$$
$$L^{(0)}(E_{2k_{2}}^{0},E_{2k_{1}}^{0};\alpha_{1},\alpha_{2})(\tau)
=(2\pi i)^{-\alpha_{1}-\alpha_{2}}
\sum_{n_{1},n_{2}>0}\frac{\sigma_{2k_{1}-1}(n_{1})\sigma_{2k_{2}-1}(n_{2})}{(n_{1}+n_{2})^{\alpha_{1}}n_{1}^{\alpha_{2}}}q^{n_{1}+n_{2}}.$$
Thus for any positive integers $i,j$ and $i+j=k$, we have the following obvious relation in $L^{2}(\mathcal{MEL_{\tau}})$:
$$L^{(0)}(E_{2k_{1}}^{0};i)(\tau)L^{(0)}(E_{2k_{2}}^{0};j)(\tau)$$
$$=\sum_{\alpha_{1}+\alpha_{2}=k}{\alpha_{1}-1 \choose i-1}L^{(0)}(E_{2k_{1}}^{0},E_{2k_{2}}^{0};\alpha_{1},\alpha_{2})(\tau)
+{\alpha_{1}-1 \choose j-1}L^{(0)}(E_{2k_{2}}^{0},E_{2k_{1}}^{0};\alpha_{1},\alpha_{2})(\tau).$$
\end{exa}

Finally, in the vector space $\mathcal{MEL_{\tau}}$, we ask $Im(\tau)>0$. If we take $\tau=0$ and $t>0$, it is obvious to see this series equals to $0$. The interesting case is $\tau=0$ and $t=0$. As stated in the introduction, it has a meromorphic extension. Denote by
$$L(E_{2k_{1}}^{0},\cdots,E_{2k_{r}}^{0};\alpha_{1},\cdots,\alpha_{r})=
L^{(0)}(E_{2k_{1}}^{0},\cdots,E_{2k_{r}}^{0};\alpha_{1},\cdots,\alpha_{r})(0)$$
which we call a multiple Eisenstein L-value. Denote by $\mathcal{MEL}$ the $\mathbb{Q}$-vector space spanned by all multiple Eisenstein L-values.

\subsection{Relations between $\mathcal{MEL_{\tau}}$ and $\mathcal{IEI_{\tau}}$}

Proposition \ref{II to MEL} below is actually a restatement of Theorem $3.2$ in \cite{Y.M1}. Since the notations and expression here are different from Manin's, we would like to restate it here.

\begin{prop}\label{II to MEL}
For any integers $r\geq1$, $k_{1},\cdots,k_{r}\geq2$ and $\alpha_{1},\cdots,\alpha_{r}\geq1$, the following formula holds:
\begin{equation*}
\begin{split}
& \ \ \ Int(E_{2k_{1}}^{0},\cdots,E_{2k_{r}}^{0};\alpha_{1},\cdots,\alpha_{r})(\tau) \\
& =\sum\limits_{i_{r}=1}^{\alpha_{r}}\sum\limits_{i_{r-1}=1}^{\alpha_{r-1}+\alpha_{r}-i_{r}}\cdots \sum\limits_{i_{1}=1}^{\alpha_{1}+\cdots+\alpha_{r}-i_{2}-\cdots-i_{r}}\frac{\alpha_{r}!}{(\alpha_{r}-i_{r}+1)!}\cdots\frac{(\alpha_{1}+\cdots+\alpha_{r}-i_{2}-\cdots-i_{r})!}{(\alpha_{1}+\cdots+\alpha_{r}-i_{1}-\cdots-i_{r}+1)!} \\
& \ \ \ \ \ \ \ \ \ \ \times(-1)^{i_{1}+\cdots+i_{r}}L^{(\alpha_{1}+\cdots+\alpha_{r}-i_{1}-\cdots-i_{r})} (E_{2k_{1}}^{0},\cdots,E_{2k_{r}}^{0};i_{1},\cdots,i_{r})(\tau).
\end{split}
\end{equation*}
\end{prop}

\noindent{\bf Proof:}
We prove the proposition by induction. Let $(1+\frac{\partial}{\partial u})^{-1}=\sum\limits_{n=0}^{\infty}(-1)^{n}\frac{\partial}{\partial u^{n}}$. Notice that for any $a,b\in \mathbb{C}$ and positive integer $\alpha$, we have
$$\int_{a}^{b}e^{u}u^{\alpha}du=e^{u}(1+\frac{\partial}{\partial u})^{-1}u^{\alpha}|_{a}^{b}.$$
Thus it is obvious to see that when $r=1$, we have:
$$Int(E_{2k_{1}}^{0};\alpha_{1})(\tau)=\sum_{i=1}^{\alpha_{1}}\frac{\alpha_{1}!}{(\alpha_{1}-i+1)!}(-1)^{i}L^{(\alpha_{1}-i)}(E_{2k_{1}}^{0};i)(\tau).$$
Now we assume that the proposition holds for $r-1$ when $r\geq2$. By induction, for $r$ we have:
\begin{equation*}
\begin{split}
&   \ \ \ \ Int(E_{2k_{1}}^{0},\cdots,E_{2k_{r}}^{0};\alpha_{1},\cdots,\alpha_{r})(\tau) \\
&    =\int_{\tau}^{i\infty}E_{2k_{1}}^{0}(\tau_{1})\tau_{1}^{\alpha_{1}-1}[\int_{\tau_{1}}^{i\infty}E_{2k_{1}}^{0}(\tau_{1})\tau_{1}^{\alpha_{1}-1}\cdots E_{2k_{r}}^{0}(\tau_{r})\tau_{r}^{\alpha_{r}-1}d\tau_{r}\cdots d\tau_{2}]d\tau_{1}\\
&    =\sum\limits_{i_{r}=1}^{\alpha_{r}}\sum\limits_{i_{r-1}=1}^{\alpha_{r-1}+\alpha_{r}-i_{r}}\cdots \sum\limits_{i_{2}=1}^{\alpha_{2}+\cdots+\alpha_{r}-i_{3}-\cdots-i_{r}}
    \frac{\alpha_{r}!}{(\alpha_{r}-i_{r}+1)!}\cdots
    \frac{(\alpha_{2}+\cdots+\alpha_{r}-i_{3}-\cdots-i_{r})!}{(\alpha_{2}+\cdots+\alpha_{r}-i_{2}-\cdots-i_{r}+1)!}\\
&    \ \ \ \ \ \ \ \ \ \ (-1)^{i_{2}+\cdots+i_{r}}\int_{\tau}^{i\infty}E_{2k_{1}}^{0}(\tau_{1})\tau_{1}^{\alpha_{1}-1}
    L^{(\alpha_{2}+\cdots+\alpha_{r}-i_{2}-\cdots-i_{r})} (E_{2k_{2}}^{0},\cdots,E_{2k_{r}}^{0};i_{2},\cdots,i_{r})(\tau_{1})
    d\tau_{1}\\
&    =\sum\limits_{i_{r}=1}^{\alpha_{r}}\sum\limits_{i_{r-1}=1}^{\alpha_{r-1}+\alpha_{r}-i_{r}}\cdots \sum\limits_{i_{1}=1}^{\alpha_{1}+\cdots+\alpha_{r}-i_{2}-\cdots-i_{r}}\frac{\alpha_{r}!}{(\alpha_{r}-i_{r}+1)!}\cdots\frac{(\alpha_{1}+\cdots+\alpha_{r}-i_{2}-\cdots-i_{r})!}{(\alpha_{1}+\cdots+\alpha_{r}-i_{1}-\cdots-i_{r}+1)!} \\
&    \ \ \ \ \ \ \ \ \ \ (-1)^{i_{1}+\cdots+i_{r}}L^{(\alpha_{1}+\cdots+\alpha_{r}-i_{1}-\cdots-i_{r})} (E_{2k_{1}}^{0},\cdots,E_{2k_{r}}^{0};i_{1},\cdots,i_{r})(\tau).
\end{split}
\end{equation*}
This completes the proof.
$\hfill\Box$\\

Proposition \ref{II to MEL} implies that we can express any iterated Eisenstein $\tau$-integral as a $\mathbb{Q}$-linear combination of multiple Eisenstein $\tau$-series, thus as $\mathbb{Q}$-vector spaces we have $\mathcal{IEI_{\tau}}\subseteq \mathcal{MEL_{\tau}}$. Conversely, by using the functor $\frac{\partial}{\partial\tau}$ and doing the induction, we have the following lemma from \cite{K.I1}:

\begin{lem}
For integers $k_{1},\cdots,k_{r}\geq2$, $\alpha_{1},\cdots,\alpha_{r}\geq1$ we have:
\begin{equation*}
\begin{split}
& \ \ \ \ (-1)^{\alpha_{1}+\cdots+\alpha_{r}}(\alpha_{1}-1)!\cdots(\alpha_{r}-1)!
L^{(0)}(E_{2k_{1}}^{0},\cdots,E_{2k_{r}}^{0};\alpha_{1},\cdots,\alpha_{r})(\tau)\\
& =\int_{\tau<z_{1}<\cdots<z_{r}<i\infty}
E_{2k_{1}}^{0}(z_{1})(z_{1}-\tau)^{\alpha_{1}-1}\cdots E_{2k_{r}}^{0}(z_{r})(z_{r}-z_{r-1})^{\alpha_{r}-1}dz_{r}\cdots dz_{1}.
\end{split}
\end{equation*}
\end{lem}

As a consequence, the following proposition holds:
\begin{prop}\label{MEL to II}
For integers $k_{1},\cdots,k_{r}\geq2$, $\alpha_{1},\cdots,\alpha_{r}\geq1$ and $t\geq0$, we have
\begin{equation*}
\begin{split}
& \ \ \ \ L^{(t)}(E_{2k_{1}}^{0},\cdots,E_{2k_{r}}^{0};\alpha_{1},\cdots,\alpha_{r})(\tau) \\
& =\frac{(-1)^{\alpha_{1}+\cdots+\alpha_{r}}}{(\alpha_{1}-1)!\cdots(\alpha_{r}-1)!}\sum\limits_{i_{1}=0}^{\alpha_{1}-1}\cdots\sum\limits_{i_{r}=0}^{\alpha_{r}-1}{\alpha_{1}-1 \choose i_{1}}\cdots{\alpha_{r}-1 \choose i_{1}}\\
& \ \ \ \ \ \ \ \ \ \ \ \ \ \ \tau^{t+i_{1}}Int(E_{2k_{1}}^{0},\cdots,E_{2k_{r}}^{0};\alpha_{1}-i_{1}+i_{2},\cdots,\alpha_{r}-i_{r}+i_{r+1})(\tau),
\end{split}
\end{equation*}
where $i_{r+1}=0$ in the above equation.
\end{prop}

If we give the length filtrations on $\mathcal{IEI_{\tau}}[\tau]$ and $\mathcal{MEL}_{\tau}$ respectively as
$$L^{n}(\mathcal{IEI_{\tau}}[\tau])=\langle \tau^{t}Int(E_{2k_{1}}^{0},\cdots,E_{2k_{r}}^{0};\alpha_{1},\cdots,\alpha_{r})(\tau) \in \mathcal{IEI_{\tau}}[\tau];r\leq n\rangle,$$
$$L^{n}(\mathcal{MEL_{\tau}})=\langle L^{(t)}(E_{2k_{1}}^{0},\cdots,E_{2k_{r}}^{0};\alpha_{1},\cdots,\alpha_{r})(\tau) \in \mathcal{MEL_{\tau}};r\leq n\rangle,$$
then combining with Proposition \ref{II to MEL} and Propostion \ref{MEL to II}, we have:

\begin{Thm}\label{bijection of vec-space}
The vector space $\mathcal{MEL_{\tau}}$ is a $\mathbb{Q}$-algebra under the normal product, and there is an algebra equation:
$$\mathcal{MEL_{\tau}}=\mathcal{IEI_{\tau}}.$$
Furthermore we have:
$$L^{n_{1}}(\mathcal{MEL_{\tau}})\cdot L^{n_{2}}(\mathcal{MEL_{\tau}})\subseteq L^{n_{1}+n_{2}}(\mathcal{MEL_{\tau}}).$$
\end{Thm}

\noindent{\bf Proof:}
We have seen that the space $\mathcal{IEI_{\tau}}[\tau]$ is a $\mathbb{Q}$-algebra. By Proposition \ref{II to MEL} we can express any element in $\mathcal{IEI_{\tau}}[\tau]$ as a $\mathbb{Q}$-linear combination of multiple Eisenstein L-series. By Proposition \ref{MEL to II} the converse statement also holds, and the expression is compatible with the filtration of them. Thus $\mathcal{MEL_{\tau}}$ is the same $\mathbb{Q}$-algebra as $\mathcal{IEI_{\tau}}$ and the last statement holds obviously.
$\hfill\Box$\\

\begin{Cor}
For integers $k_{i}\geq2$, $\alpha_{i}\geq1$ and $i=1,\cdots,r$, the iterated Eisenstein integral
$$Int(E_{2k_{1}}^{0},\cdots,E_{2k_{r}}^{0};\alpha_{1},\cdots,\alpha_{r})=Int(E_{2k_{1}}^{0},\cdots,E_{2k_{r}}^{0};\alpha_{1},\cdots,\alpha_{r})(0)$$
has a Laurent expression after the meromorphic extension.
\end{Cor}

Denote by $\mathcal{IEI}$ the $\mathbb{Q}$-vector spaces generated by all iterated Eisenstein integrals
$$\{Int(E_{2k_{1}}^{0},\cdots,E_{2k_{r}}^{0};\alpha_{1},\cdots,\alpha_{r}), k_{i}\geq2, \alpha_{i}\geq1\}.$$
There are natural reduced length filtration on $\mathcal{IEI}$ and $\mathcal{MEL}$. As a consequence of the above discussions, we have:

\begin{Cor}\label{Cor of MEL and IEI}
The vector space $\mathcal{MEL}$ and $\mathcal{IEI}$ are both $\mathbb{Q}$-algebras and we have $\mathcal{MEL}=\mathcal{IEI}$. Furthermore:
$$L^{n_{1}}(\mathcal{MEL})\cdot L^{n_{2}}(\mathcal{MEL})\subseteq L^{n_{1}+n_{2}}(\mathcal{MEL}).$$
\end{Cor}

\noindent{\bf Proof:}
The fist statement is obvious. Let $\tau=0$, $t=0$ in Proposition \ref{II to MEL} and Proposition \ref{MEL to II}, we have:
\begin{equation*}
\begin{split}
& \ \ \ Int(E_{2k_{1}}^{0},\cdots,E_{2k_{r}}^{0};\alpha_{1},\cdots,\alpha_{r}) \\
& =\sum\limits_{i_{n}=1}^{\alpha_{r}}\sum\limits_{i_{r-1}=1}^{\alpha_{r-1}+\alpha_{r}-i_{r}}\cdots \sum\limits_{i_{2}=1}^{\alpha_{1}+\cdots+\alpha_{r}-i_{3}-\cdots-i_{r}}\sum\limits_{i_{1}+\cdots+i_{r}=\alpha_{1}+\cdots+\alpha_{r}}\\
& \ \ \ \ \ \ \frac{\alpha_{r}!}{(\alpha_{r}-i_{r}+1)!}\cdots\frac{(\alpha_{1}+\cdots+\alpha_{r}-i_{2}-\cdots-i_{r})!}{(\alpha_{1}+\cdots+\alpha_{r}-i_{1}-\cdots-i_{r}+1)!}(-1)^{i_{1}+\cdots+i_{r}}L(E_{2k_{1}}^{0},\cdots,E_{2k_{r}}^{0};i_{1},\cdots,i_{r})
\end{split}
\end{equation*}
and
\begin{equation*}
\begin{split}
& \ \ \ L(E_{2k_{1}}^{0},\cdots,E_{2k_{r}}^{0};\alpha_{1},\cdots,\alpha_{r}) \\
& =\frac{1}{(\alpha_{1}-1)!\cdots(\alpha_{r}-1)!}\sum\limits_{i_{2}=0}^{\alpha_{2}-1}\cdots\sum\limits_{i_{r}=0}^{\alpha_{r}-1}{\alpha_{1}-1 \choose i_{1}}\cdots{\alpha_{r}-1 \choose i_{1}}\\
& \ \ \ \ \ \ \ \ \ \ \ \ \ \ \ \ \ \ \ \ \ Int(E_{2k_{1}}^{0},\cdots,E_{2k_{r}}^{0};\alpha_{1}-i_{1}+i_{2},\cdots,\alpha_{r}-i_{r}+i_{r+1}).
\end{split}
\end{equation*}
Thus the statement holds.
$\hfill\Box$\\

\begin{rem}
Corollary \ref{Cor of MEL and IEI} is compatible with the work of Choie and Ihara \cite{K.I1}, in where they considered cusp forms rather than Eisenstein series.
\end{rem}

At the last part of this subsection, we consider the following $\mathbb{Q}$-subspace of $\mathcal{MEL_{\tau}}$:
$$\mathcal{MEL}_{\tau}^{0}=\langle L^{(0)}(E^{0}_{2k_{1}},\cdots,E^{0}_{2k_{r}};\alpha_{1},\cdots,\alpha_{r})(\tau)\rangle.$$
We give a sketch proof that $\mathcal{MEL}_{\tau}^{0}$ is a subalgebra of $\mathcal{MEL}_{\tau}$, and thus give another way to prove that $\mathcal{MEL_{\tau}}$ is a $\mathbb{Q}$-algebra since $\mathcal{MEL_{\tau}}=\mathcal{MEL}_{\tau}^{0}[\tau]$. The fundamental tool is the formula again:
$$\frac{1}{m^{i}n^{j}}=\sum_{r+s=k}\frac{{r-1 \choose i-1}}{(m+n)^{r}n^{s}}+\frac{{r-1 \choose j-1}}{(m+n)^{r}m^{s}},$$
where $m,n\geq1$, $i,j\geq1$ are positive integers and $i+j=k$. Thus for any positive integers $n_{1},\cdots,n_{r+s}$ and $\alpha_{1},\cdots,\alpha_{r+s}$, we can express
$$\frac{1}{(n_{1}+\cdots+n_{r})^{\alpha_{1}}\cdots n_{r}^{\alpha_{r}}}\cdot\frac{1}{(n_{r+1}+\cdots+n_{r+s})^{\alpha_{r+1}}\cdots n_{r+s}^{\alpha_{s}}}$$
as a $\mathbb{Q}$-linear combination of
$$\sum_{n_{j_{1}},\cdots,n_{j_{r+s}}>0}\frac{1}{(n_{j_{1}}+\cdots+n_{j_{r+s}})^{k_{1}}\cdots n_{j_{r+s}}^{k_{r+s}}},$$
where $k_{1}+\cdots+k_{r+s}=\alpha_{1}+\cdots+\alpha_{r+s}$ and $\{j_{1},\cdots,j_{r+s}\}$ is a permutation of the set $\{1,\cdots,r+s\}$. Thus the following result holds:

\begin{prop}
The vector space $\mathcal{MEL}_{\tau}^{0}$ is $\mathbb{Q}$-algebra, and thus a subalgebra of $\mathcal{MEL}_{\tau}$.
\end{prop}

More details can be found in Lemma $1.24$ and Lemma $1.25$ in \cite{G.F}. The argument there is for multiple zeta values but still works well in this condition.

In particular, Example \ref{Stuffle relation among DEL} gives the calculation in length $2$. Notice that the relation among multiple Eisenstein $\tau$-series in this way is different from the one we get from Theorem \ref{bijection of vec-space}. They can be regarded as the analogues of stuffle relation and shuffle relation among multiple zeta values respectively.

\subsection{Linear independence}

In this subsection, we consider the linear independence of elements in $\mathcal{IEI_{\tau}}$ as complex functions of $\tau$, this is an analogue of Theorem $2.8$ in \cite{N.M}. First, we have the following theorem from \cite{M.D}:

\begin{Thm}\label{gen indepen}
Let $(\mathbf{A}, d)$ be a differential algebra over a field $k$ of characteristic $0$ with $ker(d)=k$. Assume $\mathbf{B}$ is a subfield of $\mathbf{A}$ such that $d(\mathbf{B})\subseteq \mathbf{B}$, and $X$ is any set with associated free monoid $X^{\ast}$. Suppose that $S\in \mathbf{A}\langle\langle X\rangle\rangle$ is a solution to the differential equation
$$dS=M\cdot S,$$
where $M=\sum_{x\in X}e_{x}x\in \mathbf{B}\langle\langle X\rangle\rangle$ is a homogeneous series of degree $1$, and the coefficient of the empty word in the series $S$ is supposed to be $1$. Then the following statements are equivalent:

$(1)$. The family of coefficients $\{S_{\omega}\}_{\omega\in X^{\ast}}$ of $S$ is linearly independent over $\mathbf{B}$, where $S_{\omega}$ is the coefficient of the word $\omega$ in $S$.

$(2)$. The family $\{e_{x}\}_{x\in X}$ is linearly independent over $k$, and
$$d(\mathbf{B})\cap Span_{k}\{e_{x};\ x\in X\}=\{0\}.$$
\end{Thm}

In order to simplify the expression, we may regard $E_{2k}^{0}(\tau)\tau^{\alpha-1}$ as both functions of $\tau$ and $q=e^{2\pi i\tau}$.

\begin{Thm}\label{I indepen}
The set of iterated Eisenstein $\tau$-integrals
$$\{Int(E^{0}_{2k_{1}},\cdots,E^{0}_{2k_{r}};\alpha_{1},\cdots,\alpha_{r})(\tau)\in\mathcal{IEI_{\tau}},k_{i}\geq2, \alpha_{i}\geq1, r\geq0\}$$
are $Frac(\mathbb{Z}[[e^{2\pi i\tau}]][\tau])$-linear independent functions of $\tau$.
\end{Thm}

\noindent{\bf Proof:}
We take
$$k=\mathbb{Q},\ \mathbf{A}=k[\tau]((e^{2\pi i\tau})),\ d=\frac{\partial}{\partial\tau},\ \mathbf{B}=Frac(\mathbb{Z}[[e^{2\pi i\tau}]][\tau])$$
and
$$X=\{x_{2k,\alpha}\}_{k>1,\alpha>0},\ e_{x_{2k,\alpha}}=-E_{2k}^{0}(\tau)\tau^{\alpha-1},\ M=\sum_{x_{2k,\alpha}\in X}e_{x_{2k,\alpha}}x_{2k,\alpha}.$$
These assumptions satisfy the requirement of Theorem \ref{gen indepen}. In order to prove the theorem we need to show that $\{e_{2k,\alpha}\}_{k>1,\alpha>0}$ is a set of linear independent elements and $d(\mathbf{B})\cap Span_{k}\{e_{x};\ x\in X\}=\{0\}$.

It is well-known that the set of Hecke normalized Eisenstein series $\{E_{2k}(\tau)\}$ are $\mathbb{Q}$-linear independent. Assume that
$$a_{1}E_{2k_{1}}^{0}(\tau)\tau^{\alpha_{1}-1}+\cdots+a_{n}E_{2k_{n}}^{0}(\tau)\tau^{\alpha_{n}-1}=0,\ \alpha_{i}\in \mathbb{Q}$$
where at least one of $\alpha_{i}>1$, then we have
$$a_{1}E_{2k_{1}}(\tau)\tau^{\alpha_{1}-1}+\cdots+a_{n}E_{2k_{n}}(\tau)\tau^{\alpha_{n}-1}=
a_{1}E_{2k_{1}}^{\infty}(\tau)\tau^{\alpha_{1}-1}+\cdots+a_{n}E_{2k_{n}}^{\infty}(\tau)\tau^{\alpha_{n}-1}$$
take $\tau=-\frac{1}{\tau}$, by the modular property of Eisenstein series, we have:
$$\sum_{i=1}^{n}(-1)^{\alpha_{i}-1}a_{i}E_{2k_{i}}(\tau)\tau^{2k_{i}-\alpha_{i}+1}
=\sum_{i=1}^{n}(-1)^{\alpha_{i}-1}a_{i}E_{2k_{i}}^{\infty}(\tau)\tau^{-\alpha_{i}+1}.$$

We may assume that $\alpha_{1}$ is the biggest one among $\{\alpha_{i}\}$, multiply $\tau^{\alpha_{1}-2}$ on both sides of the above equation, if $a_{1}\neq 0$, then the left-hand side function becomes a holomorphic function but the right-hand side one is not. Thus we must have $a_{1}=0$. By induction, we have $a_{i}=0$ for all $i$, thus $\{e_{2k,\alpha}\}_{k>1,\alpha>0}$ is a set of linear independent elements.

Now we use the parameter $q$ instead of $\tau$ to prove the second statement. Assume that $\sum_{i=1}^{r}a_{i}E_{2k_{i}}^{0}(q)(\frac{log(q)}{2\pi i})^{\alpha_{i}-1}\in d(\mathbf{B})$. We may furthermore assume $a_{i}\in \mathbb{Z}$, if there is a series $f\in \mathbf{B}$ such that
$$df=\sum_{i=1}^{r}a_{i}E_{2k_{i}}^{0}(q)(\frac{log(q)}{2\pi i})^{\alpha_{i}-1}\in d(\mathbf{B}).$$
By the assumption, there exists a positive integer $m$ such that $f\in \mathbb{Z}[m^{-1}]((q,log(q)))$. It means
$$\sum_{i=1}^{r}a_{i}\int_{q}^{0}E_{2k_{i}}^{0}(q)(\frac{log(q)}{2\pi i})^{\alpha_{i}-1}\frac{dq}{q}\in \mathbb{Z}[m^{-1}]((q,log(q))).$$
On the other hand, we have
$$\int_{q}^{0}E_{2k}^{0}(q)(\frac{log(q)}{2\pi i})^{\alpha_{i}-1}\frac{dq}{q}$$
$$=(2\pi i)^{-\alpha_{i}+1}
\sum_{n>0}\sigma_{2k-1}(n)\sum_{j=1}^{\alpha}(-1)^{j}\frac{(\alpha-1)\cdots(\alpha-j+2)}{n^{j}}q^{n}log(q)^{\alpha-j},$$
when $j=1$ we let $(\alpha-1)\cdots(\alpha-j+2)=1$ in the above formula. Denote by $\mathcal{G}_{2k}^{(\alpha)}=\int_{q}^{0}E_{2k}^{0}(q)(\frac{log(q)}{2\pi i})^{\alpha_{i}-1}\frac{dq}{q}$, we may further assume $\alpha_{1}=\cdots=\alpha_{r}=\alpha$ since the power of $log(q)$ in $\mathcal{G}_{2k_{i}}^{(\alpha_{i})}$ is no more than $\alpha_{i}$ and $q$ is $\mathbb{Q}(2\pi i)$-algebraic independent of $log(q)$. In this case, for any prime $p$, consider the coefficient of $q^{p}(\frac{log(q)}{2\pi i})^{\alpha-1}$ in $\mathcal{G}_{2k}^{(\alpha)}$, we have
$$\frac{\sigma_{2k-1}(p)}{p}=\frac{p^{2k-1}+1}{p}\equiv \frac{1}{p} mod(\mathbb{Z})$$
Thus for any prime $p$, $\frac{1}{p}\sum_{i=1}^{r}a_{i}\in \mathbb{Z}[m^{-1}]$, which means $\sum_{i=1}^{r}a_{i}=0$.

We may assume $k_{1}$ is the smallest number among $\{2k_{i}; i=1,\cdots,r\}$ and $a_{1}\neq 0$. Consider the coefficient of $q^{p^{2k_{1}}}(\frac{log(q)}{2\pi i})^{\alpha-1}$ in $\mathcal{G}_{2k}^{(\alpha)}$, we have
$$\frac{\sigma_{2k-1}(p^{2k_{1}})}{p^{2k-1}}=\frac{1}{p^{2k-1}}\sum_{j=0}^{2k_{1}}p^{j(2k-1)}
\equiv \{
\begin{array}{c}
\frac{1}{p^{2k_{1}}}\ mod(\mathbb{Z}) \\
\frac{1}{p^{2k_{1}}}+\frac{1}{p}\ mod(\mathbb{Z}) \\
\end{array}
$$
Thus $\frac{a_{1}}{p}+\frac{1}{p^{2k_{1}}}\sum_{i=1}^{r}a_{i}\in \mathbb{Z}[m^{-1}]$. Since $\sum_{i=1}^{r}a_{i}=0$, we have $\frac{a_{1}}{p}\in \mathbb{Z}[m^{-1}]$ for any prime $p$. This is impossible unless $a_{1}=0$, which is a contradiction.
$\hfill\Box$\\

\begin{rem}
This result is similar to Theorem $2.8$ in \cite{N.M}, but the two statements are not equivalent.
\end{rem}

As a direct consequence, we have:

\begin{Cor}
The set of elements
$$\{Int(E^{0}_{2k_{1}},\cdots,E^{0}_{2k_{r}};\alpha_{1},\cdots,\alpha_{r})(\tau)\in\mathcal{IEI_{\tau}},k_{i}\geq2, \alpha_{i}\geq1, r\geq0\}$$
are $\mathbb{C}$-linear independent functions of $\tau$.
\end{Cor}

Combining with Proposition \ref{bijection of vec-space}, by the definition of the $\mathbb{Q}$-vector spaces $\mathcal{MEL}_{\tau}$ and $\mathcal{IEI}_{\tau}$, we have:

\begin{Cor}
The following statements hold:

$(1)$. The set of elements
$$\{Int(E^{0}_{2k_{1}},\cdots,E^{0}_{2k_{r}};\alpha_{1},\cdots,\alpha_{r})(\tau)\in\mathcal{IEI_{\tau}},k_{i}\geq2, \alpha_{i}\geq1, r\geq0\}$$
form a basis of the $\mathbb{Q}$-vector space $\mathcal{IEI}_{\tau}$.

$(2)$. The set of elements
$$\{L^{(t)}(E^{0}_{2k_{1}},\cdots,E^{0}_{2k_{r}};\alpha_{1},\cdots,\alpha_{r})(\tau)\in\mathcal{MEL_{\tau}},k_{i}\geq2, \alpha_{i}\geq1,t\geq0, r\geq0\}$$
form a basis of the $\mathbb{Q}$-vector space $\mathcal{MEL}_{\tau}$.
\end{Cor}

\section{Double Eisenstein L-values}

In this section, we express a double modular value, i,e, a multiple modular value of length two (the length is naturally defined), as a $\mathbb{Q}[2\pi i]$-linear combination of double Eisenstein L-values with the help of iterated Eisenstein $\tau$-integrals.

The coefficients of $\mathcal{C}_{T}$ belong to $\mathbb{Q}[2\pi i]$ as stated in Section $2$, and thus the statement is trivial. On the other hand, directly from the definition of multiple modular values, we have
$$\mathcal{C}_{S}=\lim_{\varepsilon\rightarrow0,\eta\rightarrow i\infty}I^{\infty}(i\infty,-1/\varepsilon)|_{S}I(\varepsilon,\eta)I^{\infty}(\eta,0).$$
By induction of the length, the following proposition is straightforward.

\begin{prop}
For any given positive integers $1\leq \alpha_{i}\leq 2k_{i}-1$, $i=1,\cdots,r$, the multiple Eisenstein $L$-value
$$L(E^{0}_{2k_{1}},\cdots,E_{2k_{r}}^{0};\alpha_{1},\cdots,\alpha_{r})$$
lies in the $\mathbb{Q}$-algebra $\mathcal{MMV}[(2\pi i)^{-1}]$.
\end{prop}

In the following part of this section, we only consider the coefficients of $\mathcal{C}_{S}$. Rather then using the above expression, we would like to do this directly by definition.

\subsection{Modular property}

For $i=1,\cdots, r$ and positive integers $k_{1},\cdots ,k_{r}\geq2$, $1\leq\alpha_{i} \leq2k_{i}-1$, consider the polynomial with $2r$ variables
$$X_{1}^{2k_{1}-\alpha_{1}-1}Y_{1}^{\alpha_{1}-1}\cdots X_{r}^{2k_{r}-\alpha_{r}-1}Y_{r}^{\alpha_{r}-1},$$
which we write by $P(X,Y)$ for convenience. Denote by
$$S(2k_{1},\cdots ,2k_{r};\alpha_{1},\cdots ,\alpha_{r})$$
and
$$I(2k_{1},\cdots ,2k_{r};\alpha_{1},\cdots ,\alpha_{r})$$
the coefficients of $P(X,Y)\otimes A_{E_{2k_{1}}}\cdots A_{E_{2k_{r}}}$ in $\mathcal{C}_{S}$ and $I(i; i\infty)$ respectively. For any functions $f_{1},\cdots ,f_{r}$ and integers $\alpha_{1},\cdots ,\alpha_{r}$, denote by:
$$T(f_{1},\cdots ,f_{r}; \alpha_{1},\cdots ,\alpha_{r})=\int_{0< \tau_{1}< \cdots < \tau_{r}< i}f_{1}(\tau_{1})\tau_{1}^{\alpha_{1}}\cdots f_{r}(\tau_{r})\tau_{r}^{\alpha_{r}}\frac{d\tau_{r}}{\tau_{r}}\cdots \frac{d\tau_{1}}{\tau_{1}},$$
$$R(f_{1},\cdots ,f_{r}; \alpha_{1},\cdots ,\alpha_{r})=\int_{i< \tau_{1}< \cdots < \tau_{r}< i\infty}f_{1}(\tau_{1})\tau_{1}^{\alpha_{1}}\cdots f_{r}(\tau_{r})\tau_{r}^{\alpha_{r}}\frac{d\tau_{r}}{\tau_{r}}\cdots \frac{d\tau_{1}}{\tau_{1}}.$$
Also we will use a non-standard notation as
$$T(f; [\begin{array}{c}
                             \alpha \\
                             \beta
                           \end{array}])=\int_{0}^{i}f(\tau)(\tau^{\alpha}-\tau^{\beta})$$
and similarly for $R$ and the cases of $r>1$ to simplify the notations.

\begin{rem}
The calculations in this section are formal, they are well-defined since the regularized iterated integrals, iterated Eisenstein integrals and multiple Eisenstein L-values are well-defined (after the meromorphic extension).
\end{rem}

When $\tau\mapsto \tau+1$, it is straightforward to see
$$L^{(t)}(E^{0}_{2k_{1}},\cdots,E^{0}_{2k_{r}};\alpha_{1},\cdots,\alpha_{r})(\tau+1)=\sum_{j=0}^{t}{t \choose j} L^{(j)}(E^{0}_{2k_{1}},\cdots,E^{0}_{2k_{r}};\alpha_{1},\cdots,\alpha_{r})(\tau) $$
and
\begin{equation*}
\begin{split}
& \ \ \ Int(E_{2k_{1}}^{0},\cdots,E_{2k_{r}}^{0};\alpha_{1},\cdots,\alpha_{r})(\tau+1)\\
& =\sum\limits_{j_{1}=1}^{\alpha_{1}}\cdots\sum\limits_{j_{r}=0}^{\alpha_{r}}(-1)^{\alpha_{1}+\cdots+\alpha_{r}-j_{1}-\cdots j_{r}}{\alpha_{1} \choose \alpha_{1}-j_{1}}\cdots {\alpha_{r} \choose \alpha_{r}-j_{r}}Int(E^{0}_{2k_{1}},\cdots,E^{0}_{2k_{r}};j_{1},\cdots,j_{r})(\tau).
\end{split}
\end{equation*}

By the definition of the multiple modular values, for any $\gamma\in \Gamma$ and $\tau\in \mathfrak{h}$,
$$I(\tau; i\infty)=I(\gamma(\tau); i\infty)|_{\gamma}\mathcal{C}_{\gamma},$$
where $\mathcal{C}_{\gamma}$ is independent of the choice of $\tau$. In particular, take $\gamma=S$ and $\tau=i$, we have
$$I(i; i\infty)=I(i; i\infty)|_{\gamma}\mathcal{C}_{S}.$$
It follows that
\begin{equation*}
\begin{split}
& S(2k; \alpha)=I(2k;\alpha)-(-1)^{\alpha-1}I(2k;2k-\alpha),\\
& S(2k_{1},2k_{2}; \alpha_{1}, \alpha_{2})=I(2k_{1},2k_{2};\alpha_{1},\alpha_{2})-(-1)^{\alpha_{1}+\alpha_{1}-2}I(2k_{1},2k_{2};2k_{1}-\alpha_{1},2k_{2}-\alpha_{2})\\
& \ \ \ \ \ \ \ \ \ \ \ \ \ \ \ \ \ \ \ \ \ \ \ \ \ \ \ -(-1)^{\alpha_{1}-1}I(2k_{1};2k_{1}-\alpha_{1})S(2k_{2};\alpha_{2}).
\end{split}
\end{equation*}
By the definition of regularized iterated integral, we have:
$$I(2k;\alpha)=(-1)^{\alpha}(2\pi i)^{2k-1}{2k-2 \choose \alpha-1}[R(E_{2k}^{0};\alpha)-T(E_{2k}^{\infty};\alpha)]$$
and
\begin{equation*}
\begin{split}
& \ \ \ I(2k_{1},2k_{2};\alpha_{1} ,\alpha_{2})\\
& =(-1)^{\alpha_{1}+\alpha_{1}}(2\pi i)^{2k_{1}+2k_{2}-2}{2k_{1}-2 \choose \alpha_{1}-1}{2k_{2}-2 \choose \alpha_{2}-1}\\
& \ \ \ \ \ \ \ \ \ \ [R(E_{2k_{1}}^{0},E_{2k_{2}}^{0}; \alpha_{1},\alpha_{2})+R(E_{2k_{1}}^{\infty},E_{2k_{2}}^{0}; \alpha_{1},\alpha_{2})-R(E_{2k_{2}}^{\infty},E_{2k_{1}}^{0}; \alpha_{2},\alpha_{1})\\
& \ \ \ \ \ \ \ \ \ \ \ \ \ \ \ -R(E_{2k_{1}}^{0}; \alpha_{1})T(E_{2k_{2}}^{\infty}; \alpha_{2})+T(E_{2k_{2}}^{\infty},E_{2k_{1}}^{\infty}; \alpha_{2},\alpha_{1})].
\end{split}
\end{equation*}
The next lemma will be frequently used in the calculation:

\begin{lem}\label{fund}
With notations as above, after meromorphic extension we have:
$$R(E_{2k}^{0};\alpha)=(-1)^{\alpha}[T(E_{2k}^{0};2k-\alpha)+T(E_{2k}^{\infty};2k-\alpha)]+T(E_{2k}^{\infty};\alpha)$$
and
\begin{equation*}
\begin{split}
& \ \ \ R(E_{2k_{1}}^{\infty},E_{2k_{2}}^{0};\alpha_{1},\alpha_{2})\\
& =(-1)^{\alpha_{1}+\alpha_{2}}[T(E_{2k_{2}}^{0},E_{2k_{1}}^{\infty};2k_{2}-\alpha_{2};-\alpha_{1})+T(E_{2k_{2}}^{\infty},E_{2k_{1}}^{\infty};2k_{2}-\alpha_{2};-\alpha_{1})]\\
& \ \ \ +T(E_{2k_{2}}^{\infty},E_{2k_{1}}^{\infty};-\alpha_{2};-\alpha_{1}).
\end{split}
\end{equation*}
\end{lem}

\noindent{\bf Proof:}
We have:
$$R(E_{2k}^{0})=\int_{i}^{i\infty}(E_{2k}(\tau)-E_{2k}^{\infty}(\tau)-E_{2k}^{\infty}(\tau)\tau^{-2k}+E_{2k}^{\infty}(\tau)\tau^{-2k})\tau^{\alpha}\frac{d\tau}{\tau}.$$
Let $\tau=-\frac{1}{\tau}$, using the modular property of the Eisenstein series, we have
\begin{equation*}
\begin{split}
& R(E_{2k}^{0})=-(-1)^{\alpha-1}[\int_{0}^{i}(E_{2k}(\tau)-E_{2k}^{\infty}(\tau))\tau^{2k-\alpha}\frac{d\tau}{\tau}+\int_{0}^{i}E_{2k}^{\infty}(\tau)(\tau^{2k-\alpha}-\tau^{-\alpha})\frac{d\tau}{\tau}]\\
& \ \ \ \ \ \ \ \ \ \ =(-1)^{\alpha}[T(E_{2k}^{0};2k-\alpha)+T(E_{2k}^{\infty};2k-\alpha)-T(E_{2k}^{\infty};-\alpha)]\\
& \ \ \ \ \ \ \ \ \ \ =(-1)^{\alpha}[T(E_{2k}^{0};2k-\alpha)+T(E_{2k}^{\infty};2k-\alpha)]+T(E_{2k}^{\infty};\alpha).
\end{split}
\end{equation*}
Similarly the second formula holds.
$\hfill\Box$\\

\subsection{Calculation}

Now we calculate the following difference of two double modular values
$$S(2k_{1},2k_{2};\alpha_{1},\alpha_{2})-S(2k_{2},2k_{1};\alpha_{2},\alpha_{1}).$$
By the shuffle relation, the summation of them is easy to calculate and thus we can determine the formula of $S(2k_{1},2k_{2};\alpha_{1},\alpha_{2})$. Denote by
\begin{equation*}
\begin{split}
& A=R(E_{2k_{1}}^{0},E_{2k_{2}}^{0}; \alpha_{1},\alpha_{2})+(-1)^{\alpha_{1}+\alpha_{2}}R(E_{2k_{2}}^{0},E_{2k_{1}}^{0}; 2k_{2}-\alpha_{2},2k_{1}-\alpha_{1})\\
& \ \ \ \ \ \ -[R(E_{2k_{2}}^{0},E_{2k_{1}}^{0}; \alpha_{2},\alpha_{1})+(-1)^{\alpha_{1}+\alpha_{2}}R(E_{2k_{1}}^{0},E_{2k_{2}}^{0}; 2k_{1}-\alpha_{1},2k_{2}-\alpha_{2})],\\
& B=2[R(E_{2k_{1}}^{\infty},E_{2k_{2}}^{0}; \alpha_{1},\alpha_{2})-(-1)^{\alpha_{1}+\alpha_{2}}R(E_{2k_{1}}^{\infty},E_{2k_{2}}^{0}; 2k_{1}-\alpha_{1},2k_{2}-\alpha_{2})]\\
& \ \ \ \ \ \ -2[R(E_{2k_{2}}^{\infty},E_{2k_{1}}^{0}; \alpha_{2},\alpha_{1})-(-1)^{\alpha_{2}+\alpha_{2}}R(E_{2k_{2}}^{\infty},E_{2k_{1}}^{0}; 2k_{2}-\alpha_{2},2k_{1}-\alpha_{1})],\\
& C=-R(E_{2k_{1}}^{0};\alpha_{1})T(E_{2k_{2}}^{\infty};\alpha_{2})+(-1)^{\alpha_{1}+\alpha_{2}}R(E_{2k_{1}}^{0};2k_{1}-\alpha_{1})T(E_{2k_{2}}^{\infty};2k_{2}-\alpha_{2})\\
& \ \ \ \ \ \ +R(E_{2k_{2}}^{0};\alpha_{2})T(E_{2k_{1}}^{\infty};\alpha_{1})-(-1)^{\alpha_{1}+\alpha_{2}}R(E_{2k_{2}}^{0};2k_{2}-\alpha_{2})T(E_{2k_{1}}^{\infty};2k_{1}-\alpha_{1}),\\
& D=T(E_{2k_{2}}^{\infty},E_{2k_{1}}^{\infty};\alpha_{2},\alpha_{1})-(-1)^{\alpha_{1}+\alpha_{2}}T(E_{2k_{2}}^{\infty},E_{2k_{1}}^{\infty};2k_{2}-\alpha_{2},2k_{1}-\alpha_{1})\\
& \ \ \ \ \ \ -T(E_{2k_{1}}^{\infty},E_{2k_{2}}^{\infty};\alpha_{1},\alpha_{2})+(-1)^{\alpha_{1}+\alpha_{2}}T(E_{2k_{1}}^{\infty},E_{2k_{2}}^{\infty};2k_{1}-\alpha_{1},2k_{2}-\alpha_{2}),\\
& E=(-1)^{\alpha_{1}}I(2k_{1};2k_{1}-\alpha_{1})S(2k_{2};\alpha_{2})-(-1)^{\alpha_{2}}I(2k_{2};2k_{2}-\alpha_{2})S(2k_{1};\alpha_{1})\\
& \ \ \ =(-1)^{\alpha_{1}}I(2k_{1};2k_{1}-\alpha_{1})I(2k_{2};\alpha_{2})-(-1)^{\alpha_{2}}I(2k_{2};2k_{2}-\alpha_{2})I(2k_{1};\alpha_{1}).
\end{split}
\end{equation*}
Then by direct calculation, we have:
\begin{equation*}
\begin{split}
& \ \ \ S(2k_{1},2k_{2};\alpha_{1},\alpha_{2})-S(2k_{2},2k_{1};\alpha_{2},\alpha_{1})\\
& =(-1)^{\alpha_{1}+\alpha_{1}}(2\pi i)^{2k_{1}+2k_{2}-2}{2k_{1}-2 \choose \alpha_{1}-1}{2k_{2}-2 \choose \alpha_{2}-1}(A+B+C+D+E).
\end{split}
\end{equation*}
We will calculate it by piece. Precisely, we calculate it from higher length terms to lower ones. Note that the length means the number of the appearance of $E_{2k}^{0}(\tau)$ in the integral here since $E_{2k}^{\infty}(\tau)$ is actually a constant.

\subsubsection{Calculation of the length two part}

First, we give the following lemma to simplify the calculation.

\begin{lem}\label{Lemma of l2}
With the above notations, the following statement holds:
\begin{equation*}
\begin{split}
& \ \ \ R(E_{2k_{1}}^{0},E_{2k_{2}}^{0}; \alpha_{1},\alpha_{2})+(-1)^{\alpha_{1}+\alpha_{2}}R(E_{2k_{2}}^{0},E_{2k_{1}}^{0}; 2k_{2}-\alpha_{2},2k_{1}-\alpha_{1})\\
& =Int(E_{2k_{1}}^{0},E_{2k_{2}}^{0};\alpha_{1},\alpha_{2})+A_{2k_{1},2k_{2}}^{0}+A_{2k_{1},2k_{2}}^{'}+A_{2k_{1},2k_{2}}^{\infty},\end{split}
\end{equation*}
where
\begin{equation*}
\begin{split}
& A_{2k_{1},2k_{2}}^{0}=-T(E_{2k_{1}}^{0};\alpha_{1})R(E_{2k_{2}}^{0};\alpha_{2}),\\
& A_{2k_{1},2k_{2}}^{'}=-T(E_{2k_{1}}^{0},E_{2k_{2}}^{\infty};\alpha_{1},[\begin{array}{c}
                                                           \alpha_{2}-2k_{2}  \\
                                                           \alpha_{2}
                                                         \end{array}]
)
-T(E_{2k_{1}}^{\infty},E_{2k_{2}}^{0};[\begin{array}{c}
                                                                                                      \alpha_{1}-2k_{1} \\
                                                                                                      \alpha_{1}
                                                                                                    \end{array}
],\alpha_{2}),\\
& A_{2k_{1},2k_{2}}^{\infty}=T(E_{2k_{1}}^{\infty},E_{2k_{2}}^{\infty};[\begin{array}{c}
                                                          \alpha_{1}-2k_{1} \\
                                                          \alpha_{1}
                                                        \end{array}
],
[\begin{array}{c}
   \alpha_{2}-2k_{2} \\
   \alpha_{2}
 \end{array})
].
\end{split}
\end{equation*}
\end{lem}

\noindent{\bf Proof:}
Obviously we can divide the iterated integral as:
\begin{equation*}
\begin{split}
& \ \ \ Int(E_{2k_{1}}^{0},E_{2k_{2}}^{0};\alpha_{1},\alpha_{2}) \\
& =T(E_{2k_{1}}^{0},E_{2k_{2}}^{0};\alpha_{1},\alpha_{2})+R(E_{2k_{1}}^{0},E_{2k_{2}}^{0};\alpha_{1},\alpha_{2})+T(E_{2k_{1}};\alpha_{1})R(E_{2k_{2}};\alpha_{2}).
\end{split}
\end{equation*}
With the help of modular property of Eisenstein series, we have:
\begin{equation*}
\begin{split}
&    \ \ \ T(E_{2k_{1}}^{0},E_{2k_{2}}^{0};\alpha_{1},\alpha_{2})\\
&    =T(E_{2k_{1}},E_{2k_{2}}^{0};\alpha_{1},\alpha_{2})-T(E_{2k_{1}}^{\infty},E_{2k_{2}}^{0};\alpha_{1},\alpha_{2})\\
&    =(-1)^{\alpha_{1}+\alpha_{2}}R(E_{2k_{1}}^{0},E_{2k_{2}}^{0};2k_{1}-\alpha_{1},2k_{2}-\alpha_{2})+T(E_{2k_{1}}^{\infty},E_{2k_{2}}^{0};[\begin{array}{c}
                                                                                                                                              \alpha_{1}-2k_{1}\\
                                                                                                                                              \alpha_{1}
                                                                                                                                            \end{array}
    ],\alpha_{2})\\
&    \ \ \ +T(E_{2k_{1}}^{0},E_{2k_{2}}^{\infty};\alpha_{1},\alpha_{2}-2k_{2})-T(E_{2k_{1}}^{0},E_{2k_{2}}^{\infty};\alpha_{1},\alpha_{2})\\
&    \ \ \ -T(E_{2k_{1}}^{\infty},E_{2k_{2}}^{\infty};[\begin{array}{c}
                                                         \alpha_{1}-2k_{1} \\
                                                         \alpha_{1}
                                                       \end{array}
    ],[\begin{array}{c}
         \alpha_{2}-2k_{2} \\
         \alpha_{2}
       \end{array}
    ])\\
&    =(-1)^{\alpha_{1}+\alpha_{2}}R(E_{2k_{1}}^{0},E_{2k_{2}}^{0};2k_{1}-\alpha_{1},2k_{2}-\alpha_{2})-A_{2k_{1},2k_{2}}^{'}-A_{2k_{1},2k_{2}}^{\infty}.
\end{split}
\end{equation*}
This implies
\begin{equation*}
\begin{split}
& \ \ \ R(E_{2k_{1}}^{0},E_{2k_{2}}^{0}; \alpha_{1},\alpha_{2})+(-1)^{\alpha_{1}+\alpha_{2}}R(E_{2k_{2}}^{0},E_{2k_{1}}^{0}; 2k_{2}-\alpha_{2},2k_{1}-\alpha_{1})\\
& =Int(E_{2k_{1}}^{0},E_{2k_{2}}^{0};\alpha_{1};\alpha_{2})+A_{2k_{1},2k_{2}}^{0}+A_{2k_{1},2k_{2}}^{'}+A_{2k_{1},2k_{2}}^{\infty}
\end{split}
\end{equation*}
and thus the lemma holds.
$\hfill\Box$\\

According to Lemma \ref{Lemma of l2}, we can write
$$A=Int(E_{2k_{1}}^{0},E_{2k_{2}}^{0};\alpha_{1},\alpha_{2})-Int(E_{2k_{2}}^{0},E_{2k_{1}}^{0};\alpha_{2},\alpha_{1})+A^{0}+A'+A^{\infty},$$ where $A^{*}=A_{2k_{1},2k_{2}}^{*}-A_{2k_{2},2k_{1}}^{*}$ for $*\in\{0,', \infty\}$. Then the question reduces to the calculation of the sum
$$A^{0}+A'+A^{\infty}+B+C+D+E.$$
Then according to Lemma \ref{fund}, we can rewrite $A^{0}$ and $E$ as
\begin{equation*}
\begin{split}
& A^{0}=-(-1)^{\alpha_{2}}T(E_{2k_{1}}^{0};\alpha_{1})[T(E_{2k_{2}}^{0};2k_{2}-\alpha_{2})+T(E_{2k_{2}}^{\infty};[\begin{array}{c}
                                                                                                                    2k_{2}-\alpha_{2} \\
                                                                                                                    \alpha_{2}
                                                                                                                  \end{array}
])]\\
& \ \ \ \ \ \ \ \ +(-1)^{\alpha_{1}}T(E_{2k_{2}}^{0};\alpha_{2})[T(E_{2k_{1}}^{0};2k_{1}-\alpha_{1})+T(E_{2k_{1}}^{\infty};[\begin{array}{c}
                                                                                                                    2k_{1}-\alpha_{1} \\
                                                                                                                    \alpha_{1}
                                                                                                                  \end{array}
])],\\
& E=(-1)^{\alpha_{2}}[T(E_{2k_{1}}^{0};\alpha_{1})+T(E_{2k_{1}}^{\infty};\alpha_{1})][T(E_{2k_{2}}^{0};2k_{2}-\alpha_{2})+T(E_{2k_{2}}^{\infty};2k_{2}-\alpha_{2})]\\
& \ \ \ \ \ \ \ -(-1)^{\alpha_{1}}[T(E_{2k_{2}}^{0};\alpha_{2})+T(E_{2k_{2}}^{\infty};\alpha_{2})][T(E_{2k_{1}}^{0};2k_{1}-\alpha_{1})+T(E_{2k_{1}}^{\infty};2k_{1}-\alpha_{1})].
\end{split}
\end{equation*}
Then we have
\begin{equation*}
\begin{split}
&    \ \ \ A^{0}+E \\
&    =-T(E_{2k_{1}}^{0};\alpha_{1})T(E_{2k_{2}}^{\infty};\alpha_{2})+T(E_{2k_{2}}^{0};\alpha_{2})T(E_{2k_{1}}^{\infty};\alpha_{1}) \\
&    \ \ \ -(-1)^{\alpha_{1}}T(E_{2k_{1}}^{0};2k_{1}-\alpha_{1})T(E_{2k_{2}}^{\infty};\alpha_{2})+(-1)^{\alpha_{2}}T(E_{2k_{2}}^{0};2k_{2}-\alpha_{2})T(E_{2k_{1}}^{\infty};\alpha_{1})\\
&    \ \ \ +(A^{0}+E)^{\infty}
\end{split}
\end{equation*}
where
\begin{equation*}
\begin{split}
& \ \ \ (A^{0}+E)^{\infty}\\
& =-(-1)^{\alpha_{1}}T(E_{2k_{1}}^{\infty};2k_{1}-\alpha_{1})T(E_{2k_{2}}^{\infty};\alpha_{2})+(-1)^{\alpha_{2}}T(E_{2k_{2}}^{\infty};2k_{2}-\alpha_{2})T(E_{2k_{1}}^{\infty};\alpha_{1}).
\end{split}
\end{equation*}
Up to now, we have determined the length $2$ part of the difference
$$S(2k_{1},2k_{2};\alpha_{1},\alpha_{2})-S(2k_{2},2k_{1};\alpha_{2},\alpha_{1}).$$

\subsubsection{Calculation of length one part}

Use Lemma \ref{fund} again, we can rewrite
\begin{equation*}
\begin{split}
& B=2(-1)^{\alpha_{1}+\alpha_{2}}[T(E_{2k_{2}}^{0},E_{2k_{1}}^{\infty};2k_{2}-\alpha_{2},-\alpha_{1})+T(E_{2k_{2}}^{\infty},E_{2k_{1}}^{\infty};[\begin{array}{c}
                                                                                                                                                   2k_{2}-\alpha_{2} \\
                                                                                                                                                   -\alpha_{2}
                                                                                                                                                 \end{array}
],-\alpha_{1})\\
& \ \ \ \ \ \ \ -T(E_{2k_{1}}^{0},E_{2k_{2}}^{\infty};2k_{1}-\alpha_{1},-\alpha_{2})-T(E_{2k_{1}}^{\infty},E_{2k_{2}}^{\infty};[\begin{array}{c}
                                                                                                                                                   2k_{1}-\alpha_{1} \\
                                                                                                                                                   -\alpha_{1}
                                                                                                                                                 \end{array}
],-\alpha_{2})]\\
& \ \ \ \ \ \ \ -2[T(E_{2k_{2}}^{0},E_{2k_{1}}^{\infty};\alpha_{2},\alpha_{1}-2k_{1})+T(E_{2k_{2}}^{\infty},E_{2k_{1}}^{\infty};[\begin{array}{c}
                                                                                                                                                   \alpha_{2} \\
                                                                                                                                                   \alpha_{2}-2k_{2}
                                                                                                                                                 \end{array}
],\alpha_{1}-2k_{1})\\
& \ \ \ \ \ \ \ -T(E_{2k_{1}}^{0},E_{2k_{2}}^{\infty};\alpha_{1},\alpha_{2}-2k_{2})-T(E_{2k_{1}}^{\infty},E_{2k_{2}}^{\infty};[\begin{array}{c}
                                                                                                                                                   \alpha_{1} \\
                                                                                                                                                   \alpha_{1}-2k_{1}
                                                                                                                                                 \end{array}
],\alpha_{2}-2k_{2})],
\end{split}
\end{equation*}
and
\begin{equation*}
\begin{split}
& A'=-T(E_{2k_{1}}^{\infty};[\begin{array}{c}
                              \alpha_{1}-2k_{1} \\
                              \alpha_{1}
                            \end{array}]
)T(E_{2k_{2}}^{0};\alpha_{2})
+T(E_{2k_{2}}^{\infty};[\begin{array}{c}
                              \alpha_{2}-2k_{2} \\
                              \alpha_{2}
                            \end{array}]
)T(E_{2k_{1}}^{0};\alpha_{1})\\
& \ \ \ \ \ \ \ -2T(E_{2k_{1}}^{0},E_{2k_{2}}^{\infty};\alpha_{1},[\begin{array}{c}
                                                       \alpha_{2}-2k_{2} \\
                                                       \alpha_{2}
                                                     \end{array}
])
+2T(E_{2k_{2}}^{0},E_{2k_{1}}^{\infty};\alpha_{2},[\begin{array}{c}
                                                       \alpha_{1}-2k_{1} \\
                                                       \alpha_{1}
                                                     \end{array}
]).
\end{split}
\end{equation*}
Thus we have:
\begin{equation*}
\begin{split}
&    \ \ \ A^{'}+B \\
&    =-T(E_{2k_{1}}^{\infty};[\begin{array}{c}
                              \alpha_{1}-2k_{1} \\
                              \alpha_{1}
                            \end{array}]
)T(E_{2k_{2}}^{0};\alpha_{2})
+T(E_{2k_{2}}^{\infty};[\begin{array}{c}
                              \alpha_{2}-2k_{2} \\
                              \alpha_{2}
                            \end{array}]
)T(E_{2k_{1}}^{0};\alpha_{1}) \\
&    \ \ \ +2\eta(2k_{1},2k_{2};\alpha_{1},\alpha_{2})-2\eta(2k_{2},2k_{1};\alpha_{2},\alpha_{1})+B^{\infty},
\end{split}
\end{equation*}
where
$$\eta(2k_{1},2k_{2};\alpha_{1},\alpha_{2})=T(E_{2k_{1}}^{0},E_{2k_{2}}^{\infty};\alpha_{1},\alpha_{2})-(-1)^{\alpha_{1}+\alpha_{2}}T(E_{2k_{1}}^{0},E_{2k_{2}}^{\infty};2k_{1}-\alpha_{1},-\alpha_{2})$$
and $B^{\infty}$ is the rest part of $B$ of lower length as:
\begin{equation*}
\begin{split}
& B^{\infty}=2(-1)^{\alpha_{1}+\alpha_{2}}[T(E_{2k_{2}}^{\infty},E_{2k_{1}}^{\infty};[\begin{array}{c}
                                                                                                                                                   2k_{2}-\alpha_{2} \\
                                                                                                                                                   -\alpha_{2}
                                                                                                                                                 \end{array}
],-\alpha_{1})-T(E_{2k_{1}}^{\infty},E_{2k_{2}}^{\infty};[\begin{array}{c}
                                                                                                                                                   2k_{1}-\alpha_{1} \\
                                                                                                                                                   -\alpha_{1}
                                                                                                                                                 \end{array}
],-\alpha_{2})]\\
& \ \ \ \ \ \ \ \ -2[T(E_{2k_{2}}^{\infty},E_{2k_{1}}^{\infty};[\begin{array}{c}
                                                                                                                                                   \alpha_{2} \\
                                                                                                                                                   \alpha_{2}-2k_{2}
                                                                                                                                                 \end{array}
],\alpha_{1}-2k_{1})-T(E_{2k_{1}}^{\infty},E_{2k_{2}}^{\infty};[\begin{array}{c}
                                                                                                                                                   \alpha_{1} \\
                                                                                                                                                   \alpha_{1}-2k_{1}
                                                                                                                                                 \end{array}
],\alpha_{2}-2k_{2})].
\end{split}
\end{equation*}

\begin{lem}\label{yt}
We have the following formula:
\begin{equation*}
\begin{split}
& \eta(2k_{1},2k_{2};\alpha_{1},\alpha_{2})=T(E_{2k_{2}}^{\infty};\alpha_{2})[T(E_{2k_{1}}^{0};\alpha_{1})+(-1)^{\alpha_{1}}T(E_{2k_{1}}^{0};2k_{1}-\alpha_{1})]\\
& \ \ \ \ \ \ \ \ \ \ \ \ \ \ \ \ \ \ \ \ \ \ \ \ -\frac{E_{2k_{2}}^{\infty}(\tau_{2})}{\alpha_{2}}[T(E_{2k_{1}}^{0};\alpha_{1}+\alpha_{2})+(-1)^{\alpha_{1}+\alpha_{2}}T(E_{2k_{1}}^{0};2k_{1}-\alpha_{1}-\alpha_{2})].
\end{split}
\end{equation*}
\end{lem}

\noindent{\bf Proof:}
We just need to notice the following two facts:
\begin{equation*}
\begin{split}
& T(E_{2k}^{\infty};\alpha)=E_{2k}^{\infty}(\tau)\frac{i^{\alpha}}{\alpha},\\
& T(E_{2k_{1}}^{0},E_{2k_{2}}^{\infty};\alpha_{1},\alpha_{2})
=\frac{E_{2k_{2}}^{\infty}(\tau_{2})}{\alpha_{2}}[i^{\alpha_{2}}T(E_{2k_{1}}^{0};\alpha_{1})-T(E_{2k_{1}}^{0};\alpha_{1}+\alpha_{2})],
\end{split}
\end{equation*}
and then the lemma holds from direct calculation.
$\hfill\Box$\\

On the other hand,
\begin{equation*}
\begin{split}
& \ \ \ C-C^{\infty}\\
& =-(-1)^{\alpha_{1}}T(E_{2k_{1}}^{0};2k_{1}-\alpha_{1})T(E_{2k_{2}}^{\infty};\alpha_{2})+(-1)^{\alpha_{2}}T(E_{2k_{2}}^{0};2k_{2}-\alpha_{2})T(E_{2k_{1}}^{\infty};\alpha_{1})\\
& \ \ \ +(-1)^{\alpha_{2}}T(E_{2k_{1}}^{0};\alpha_{1})T(E_{2k_{2}}^{\infty};2k_{2}-\alpha_{2})-(-1)^{\alpha_{1}}T(E_{2k_{2}}^{0};\alpha_{2})T(E_{2k_{1}}^{\infty};2k_{1}-\alpha_{1}),
\end{split}
\end{equation*}
where
\begin{equation*}
\begin{split}
& C^{\infty}=-(-1)^{\alpha_{1}}T(E_{2k_{1}}^{\infty};[\begin{array}{c}
                                                        2k_{1}-\alpha_{1} \\
                                                        -\alpha_{1}
                                                      \end{array}
])T(E_{2k_{2}}^{\infty};\alpha_{2})+(-1)^{\alpha_{2}}T(E_{2k_{1}}^{\infty};[\begin{array}{c}
                                                        \alpha_{1} \\
                                                        2k_{1}-\alpha_{1}
                                                      \end{array}
])T(E_{2k_{2}}^{\infty};2k_{2}-\alpha_{2})\\
& \ \ \ \ \ \ \ \ +(-1)^{\alpha_{2}}T(E_{2k_{2}}^{\infty};[\begin{array}{c}
                                                                                               2k_{2}-\alpha_{2} \\
                                                                                               -\alpha_{2}
                                                                                             \end{array}]
)T(E_{2k_{1}}^{\infty};\alpha_{1})-(-1)^{\alpha_{1}}T(E_{2k_{2}}^{\infty};[\begin{array}{c}
                                                                                               \alpha_{2} \\
                                                                                               2k_{2}-\alpha_{2}
                                                                                             \end{array}]
)T(E_{2k_{1}}^{\infty};2k_{1}-\alpha_{1}).
\end{split}
\end{equation*}
Now let us sum all the above terms up, combining with Lemma \ref{yt}, we have:

\begin{equation*}
\begin{split}
&    \ \ \ A^{0}+A^{'}+B+C+E \\
&    =-2[T(E_{2k_{1}}^{0};\alpha_{1})+(-1)^{\alpha_{1}}T(E_{2k_{1}}^{0};2k_{1}-\alpha_{1})]T(E_{2k_{2}}^{\infty};\alpha_{2}) \\
&    \ \ \ +2[T(E_{2k_{2}}^{0};\alpha_{2})+(-1)^{\alpha_{2}}T(E_{2k_{2}}^{0};2k_{2}-\alpha_{2})]T(E_{2k_{1}}^{\infty};\alpha_{1}) \\
&    \ \ \ +2\eta(2k_{1},2k_{2};\alpha_{1},\alpha_{2})-2\eta(2k_{2},2k_{1};\alpha_{2},\alpha_{1})+(A^{0}+E)^{\infty}+B^{\infty}+C^{\infty}\\
&    =-2\frac{E_{2k_{2}}^{\infty}(\tau_{2})}{\alpha_{2}}[T(E_{2k_{1}}^{0};\alpha_{1}+\alpha_{2})+(-1)^{\alpha_{1}+\alpha_{2}}T(E_{2k_{1}}^{0};2k_{1}-\alpha_{1}-\alpha_{2})]\\
&    \ \ \ +2\frac{E_{2k_{1}}^{\infty}(\tau_{1})}{\alpha_{1}}[T(E_{2k_{2}}^{0};\alpha_{1}+\alpha_{2})+(-1)^{\alpha_{1}+\alpha_{2}}T(E_{2k_{2}}^{0};2k_{2}-\alpha_{1}-\alpha_{2})]\\
&    \ \ \ +(A^{0}+E)^{\infty}+B^{\infty}+C^{\infty}\\
&    =-2\frac{E_{2k_{2}}^{\infty}(\tau_{2})}{\alpha_{2}}[Int(E_{2k_{1}}^{0};\alpha_{1}+\alpha_{2})+(-1)^{\alpha_{1}+\alpha_{2}}T(E_{2k_{1}}^{\infty};[\begin{array}{c}
                                                                                                              -\alpha_{1}-\alpha_{2} \\
                                                                                                              2k_{1}-\alpha_{1}-\alpha_{2}
                                                                                                            \end{array}
    ])]\\
&    \ \ \ +2\frac{E_{2k_{1}}^{\infty}(\tau_{1})}{\alpha_{1}}[Int(E_{2k_{2}}^{0};\alpha_{1}+\alpha_{2})+(-1)^{\alpha_{1}+\alpha_{2}}T(E_{2k_{2}}^{\infty};[\begin{array}{c}
                                                                                                              -\alpha_{1}-\alpha_{2} \\
                                                                                                              2k_{2}-\alpha_{1}-\alpha_{2}
                                                                                                            \end{array}
    ])]\\
&    \ \ \ +(A^{0}+E)^{\infty}+B^{\infty}+C^{\infty}.
\end{split}
\end{equation*}
Thus we have determined the length two and length one parts of the difference. In order to calculate the difference, we only need to calculate the rest part of the sum as
\begin{equation*}
\begin{split}
&    \ \ \ (A^{0}+E)^{\infty}+A^{\infty}+B^{\infty}+C^{\infty}+D \\
&    =2T(E_{2k_{2}}^{\infty},E_{2k_{1}}^{\infty};\alpha_{2},\alpha_{1}-2k_{1})-2T(E_{2k_{1}}^{\infty},E_{2k_{2}}^{\infty};\alpha_{1},\alpha_{2}-2k_{2})\\
&    \ \ \ +2T(E_{2k_{1}}^{\infty},E_{2k_{2}}^{\infty};\alpha_{1},\alpha_{2})-2T(E_{2k_{2}}^{\infty},E_{2k_{1}}^{\infty};\alpha_{2},\alpha_{1}).\\
\end{split}
\end{equation*}
This is a constant.

\begin{rem}
Note that $T(E_{2k_{1}}^{\infty},E_{2k_{2}}^{\infty};\alpha_{1},\alpha_{2})=T(E_{2k_{2}}^{\infty},E_{2k_{1}}^{\infty};\alpha_{1},\alpha_{2})$, this will make the calculation much easier.
\end{rem}

\subsubsection{Main results}

As a summary, we can state the main result in this section.

\begin{Thm}
With the above notations, we have
\begin{equation*}
\begin{split}
&    \ \ \ S(2k_{1},2k_{2};\alpha_{1},\alpha_{2})-S(2k_{2},2k_{1};\alpha_{2},\alpha_{1}) \\
&    =(-1)^{\alpha_{1}+\alpha_{2}}(2\pi i)^{2k_{1}+2k_{2}-2}{2k_{1}-2 \choose \alpha_{1}-1}{2k_{2}-2 \choose \alpha_{2}-1}\\
&    \ \ \ \ \ \ \{Int(E_{2k_{1}}^{0},E_{2k_{2}}^{0};\alpha_{1},\alpha_{2})-\frac{b_{2k_{2}}}{2k_{2}\alpha_{2}}Int(E_{2k_{1}}^{0};\alpha_{1}+\alpha_{2})\\
&    \ \ \ \ \ \ \ \ \ -[Int(E_{2k_{2}}^{0},E_{2k_{1}}^{0};\alpha_{2},\alpha_{1})-\frac{b_{2k_{1}}}{2k_{1}\alpha_{1}}Int(E_{2k_{2}}^{0};\alpha_{1}+\alpha_{2})]+\frac{b_{2k_{1}}b_{2k_{2}}(\alpha_{2}-\alpha_{1})}{8k_{1}k_{2}\alpha_{1}\alpha_{2}(\alpha_{1}+\alpha_{2})}\}.
\end{split}
\end{equation*}
\end{Thm}

\noindent{\bf Proof:}
We just need to notice that
$$T(E_{2k_{1}}^{\infty},E_{2k_{2}}^{\infty};\alpha_{1},\alpha_{2})=\frac{i^{\alpha_{1}\alpha_{2}}E_{2k_{1}}^{\infty}(\tau_{1})E_{2k_{2}}^{\infty}(\tau_{2})}{\alpha_{1}(\alpha_{1}+\alpha_{2})}$$
and that $E_{2k}^{\infty}(\tau)=-\frac{b_{2k}}{4k}$, where $b_{2k}$ is the $2k$-th Bernoulli number. Then the theorem follows from the above discussion.
$\hfill\Box$\\

Since $S(2k;\alpha)=-(-1)^{\alpha}(2\pi i)^{2k-1}{2k-2 \choose \alpha-1}Int(E_{2k}^{0};\alpha)$, denote by
$$C_{\alpha_{1},\alpha_{2}}^{k_{1},k_{2}}= (-1)^{\alpha_{1}+\alpha_{2}}(2\pi i)^{2k_{1}+2k_{2}-2}{2k_{1}-1 \choose \alpha_{1}-1}{2k_{2}-1 \choose \alpha_{2}-1},$$
by the shuffle product of iterated integrals, we have
\begin{equation*}
\begin{split}
& \ \ \ S(2k_{1},2k_{2};\alpha_{1},\alpha_{2})+S(2k_{2},2k_{1};\alpha_{2},\alpha_{1})\\
& =S(2k_{1};\alpha_{1})S(2k_{2};\alpha_{2})\\
& =C_{\alpha_{1},\alpha_{2}}^{k_{1},k_{2}}\times[Int(E_{2k_{1}}^{0},E_{2k_{2}}^{0};\alpha_{1},\alpha_{2})+ Int(E_{2k_{1}}^{0},E_{2k_{2}}^{0};\alpha_{1},\alpha_{2})].
\end{split}
\end{equation*}
Combining with Proposition \ref{II to MEL}, we give the double modular value a linear combination of the double Eisenstein L-values as below:

\begin{Cor}\label{DMV to DEL}
The following formula holds:
\begin{equation*}
\begin{split}
& \ \ \ S(2k_{1},2k_{2};\alpha_{1},\alpha_{2})\\
& = C_{\alpha_{1},\alpha_{2}}^{k_{1},k_{2}}\times[Int(E_{2k_{1}}^{0},E_{2k_{2}}^{0};\alpha_{1},\alpha_{2})\\
& \ \ \ \ \ \ -\frac{b_{2k_{2}}}{4k_{2}\alpha_{2}}Int(E^{0}_{2k_{1}};\alpha_{1}+\alpha_{2})+ \frac{b_{2k_{1}}}{4k_{1}\alpha_{1}}Int(E^{0}_{2k_{2}};\alpha_{1}+\alpha_{2})+\frac{b_{2k_{1}}b_{2k_{2}}(\alpha_{2}-\alpha_{1})}{16k_{1}k_{2}\alpha_{1}\alpha_{2}(\alpha_{1}+\alpha_{2})}]\\
& =(-1)^{\alpha_{1}+\alpha_{2}} C_{\alpha_{1},\alpha_{2}}^{k_{1},k_{2}}\times[\sum\limits_{i_{1}+i_{2}=\alpha_{1}+\alpha_{2};1\leq i_{2}\leq \alpha_{2}}\frac{\alpha_{2}!i_{1}!}{(\alpha_{2}-i_{2}+1)!}L(E_{2k_{1}}^{0},E_{2k_{2}}^{0};i_{1},i_{2})\\
& \ \ \ \ \ \ -\frac{b_{2k_{2}}(\alpha_{1}+\alpha_{2})!}{4k_{2}\alpha_{2}}L(E^{0}_{2k_{1}};\alpha_{1}+\alpha_{2})+ \frac{b_{2k_{1}}(\alpha_{1}+\alpha_{2})!}{4k_{1}\alpha_{1}}L(E^{0}_{2k_{2}};\alpha_{1}+\alpha_{2})+\frac{b_{2k_{1}}b_{2k_{2}}(\alpha_{2}-\alpha_{1})}{16k_{1}k_{2}\alpha_{1}\alpha_{2}(\alpha_{1}+\alpha_{2})}].
\end{split}
\end{equation*}
\end{Cor}

\subsection{Some remarks}

In the final part of this section, we give some remarks about the double Eisenstein L-values.

First, by corollary \ref{DMV to DEL}, one finds that the Riemann zeta value $\zeta(\alpha_{1}+\alpha_{2})$ appears in $S(2k_{1},2k_{2};\alpha_{1},\alpha_{2})$.

Next, with the notations in Theorem \ref{Brown's l2 calculation}, for any $k_{1},k_{2}\geq2$ and $k=0$, we have $I_{2k_{1},2k_{2}}^{0}(S)=\pi_{d}(Im(P))$ for
\begin{equation*}
\begin{split}
& P=\sum\limits_{\alpha_{1}=1}^{2k_{1}-1}\sum\limits_{\alpha_{2}=1}^{2k_{2}-1}S(2k_{1},2k_{2};\alpha_{1},\alpha_{2})X_{1}^{2k_{1}-\alpha_{1}-1}Y_{1}^{\alpha_{1}-1} X_{2}^{2k_{2}-\alpha_{2}-1}Y_{2}^{\alpha_{2}-1}\\
& \ \ \ \ \ \ +\frac{(2k_{1}-2)!(2k_{2}-2)!}{4}(2\pi i)^{2k_{2}-1}\zeta(2k_{1}-1)\sum\limits_{i=1}^{k_{2}-1}\frac{b_{2i}b_{2k_{2}-2i}}{(2i)!(2k_{2}-2i)!}X_{1}^{2k_{1}-2}X_{2}^{2i-1}Y_{2}^{2k_{2}-2i-1}\\
& \ \ \ \ \ \ +\frac{(2k_{1}-2)!(2k_{2}-2)!}{4}(2\pi i)^{2k_{1}-1}\zeta(2k_{2}-1)\sum\limits_{j=1}^{k_{1}-1}\frac{b_{2j}b_{2k_{1}-2j}}{(2j)!(2k_{1}-2j)!}X_{1}^{2j-1}Y_{1}^{2k_{1}-2j-1}Y_{2}^{2k_{2}-2},
\end{split}
\end{equation*}
where $Im(P)$ means the imaginary part of $P$ and
$$\pi_{d}:\mathbb{Q}[X_{1},X_{2},Y_{1},Y_{2}]\rightarrow\mathbb{Q}[X,Y];X_{i}\mapsto X,Y_{i}\mapsto Y,i=1,2.$$
On the other hand, keep the notations as Theorem \ref{Brown's l2 calculation}, we have
$$I^{0}_{2a,2b}(S)\equiv \sum_{\{g\}}\Lambda(g,w+k)P_{g}^{-}(X,Y)\ (mod\ \delta^{0}(V_{w-2}\otimes \mathbb{C})_{S}).$$
Actually, by Theorem $9.4$ in \cite{F.B2}, this is an equation rather then a modulo equivalence. By comparing the coefficients of the polynomials, one finds that the multiple Eisenstein L-value $L(E_{2k_{1}}^{0},E_{2k_{2}}^{0};\alpha_{1},\alpha_{2})$ is related to the classical L-values of cusp forms of weight $2k_{1}+2k_{2}-2$ outside the critical line by Theorem \ref{Brown's l2 calculation}.

In particular, when $\alpha_{1}=\alpha_{2}=1$, denote by $w=2k_{1}+2k_{2}-2$, we have
$$(2\pi i)^{w}(L(E_{2k_{1}}^{0},E_{2k_{2}}^{0};1,1))-\frac{b_{2k_{2}}}{2k_{2}}L(E_{2k_{1}}^{0};2)+ \frac{b_{2k_{1}}}{2k_{1}}L(E_{2k_{2}}^{0};2)=\delta_{k_{1},k_{2}}$$
for some constant number $\delta_{k_{1},k_{2}}\in\mathbb{C}$. We can see that $L(E_{2k}^{0};2)\in\mathbb{Q}[2\pi i]$. As suggested by Brown \cite{F.B2}, the constant $\delta_{k_{1},k_{2}}$ should be related to the double zeta value $\zeta(2k_{1}-1,2k_{2}-1)$.

Finally, combining with Example \ref{Stuffle relation among DEL} and Corollary \ref{DMV to DEL}, one may get another relationship among double modular values rather then the one coming from shuffle property of (regularized) iterated integrals.

\section*{Acknowledgment}

I express my sincere gratitude to Matthes for pointing out some errors in my previous manuscript. Also, I appreciate Matthes and Choie for giving the excellent suggestions, which help a lot to improve this paper.

\end{document}